\documentclass[a4paper,reqno]{amsart}

\usepackage[utf8]{inputenc}
\usepackage{amssymb}
\usepackage{enumitem}
\usepackage{mathrsfs}
\usepackage{mathtools}
\usepackage{amsmath}
\usepackage{xcolor}
\usepackage[colorlinks=true]{hyperref}
\hypersetup{linkcolor=violet,urlcolor=cyan,citecolor=cyan}

\setlist[enumerate]{label=\emph{(\roman*)}}

\usepackage{environ}

\newtheorem{theorem}{Theorem}[section]
\newtheorem{corollary}[theorem]{Corollary}
\newtheorem{lemma}[theorem]{Lemma}
\newtheorem{proposition}[theorem]{Proposition}

\theoremstyle{definition}

\newtheorem{remark}[theorem]{Remark}
\numberwithin{equation}{section}
\newcommand{\R}{\mathbb{R}}

\def \C {{\mathbb{C}}}
\def \d {{\rm{d}}}
\def \pt{\partial_{t}}

\begin{document}

\parindent=0pt

	\title[The 3D DKG system with uniform energy]
	{Global solution to the 3D Dirac--Klein-Gordon system with uniform energy bounds}
	
		\author[S.~Dong]{Shijie Dong}
	\address{Southern University of Science and Technology, SUSTech International Center for Mathematics, and Department of Mathematics, 518055 Shenzhen, P.R. China.}
\email{dongsj@sustech.edu.cn, shijiedong1991@hotmail.com}
	
		\author[K.~Li]{Kuijie Li}
\address{Nankai University, School of Mathematical Sciences and LPMC, Tianjin, 300071, P.R. China.}
\email{kuijiel@nankai.edu.cn}

	\author[X.~Yuan]{Xu Yuan}
	\address{Department of Mathematics, The Chinese University of Hong Kong, Shatin, N.T., Hong Kong, P.R. China.}
	\email{xu.yuan@cuhk.edu.hk}
	\begin{abstract}
		On the (1+3) dimensional Minkowski spacetime, for small, regular initial data, it is well-known that the Dirac-Klein-Gordon system admits a global solution. In the present paper, we aim to establish the uniform boundedness of the total energy of the solution for this system. The proof relies on Klainerman's vector field and Alinhac's ghost weight methods. 
		
		The main difficulty originates from the slow decay nature of the Dirac and wave components in three space dimensions. To overcome the difficulty, a sharp understanding of the structure for this system, and a new weighted conformal energy estimate are required. In addition, we also provide a few scattering results for the system.
	\end{abstract}
	\maketitle
\section{Introduction}
\subsection{The model problems}
The 3D Dirac--Klein-Gordon system arises as a basic model in relativistic quantum mechanics, describing a scalar or pseudoscalar field and a Dirac field through Yukawa interactions.
The model of interest can be expressed as
\begin{equation}\label{eq:DKG-Mm}
\left\{\begin{aligned}
&-i\gamma^{\mu}\partial_{\mu}\psi + M \psi =\phi V\psi,
\\
&-\Box\phi + m^2 \phi =\psi^{*}U\psi.
\end{aligned}\right.
\end{equation}
The unknown $\psi$ takes values in $\mathbb{C}^4$, and the unknown $\phi$ takes values in $\mathbb{R}$. We use $M, m$ to denote the masses of the Dirac field $\psi$ and the Klein-Gordon field $\phi$, respectively. The Greek indices $\left\{\mu,\nu,\dots\right\}$ range over $\left\{0,1,2,3\right\}$, Roman indices $\left\{a,b,\dots\right\}$ range over $\left\{1,2,3\right\}$, and the Einstein summation convention for repeated upper and lower indices is adopted.
In the above, $V$ and $U$ are $4\times 4$ matrices, and the other $4\times 4$ matrices $\left\{\gamma^{0},\gamma^{1},\gamma^{2},\gamma^{3}\right\}$ are Dirac matrices, which verify the relations:
\begin{equation}\label{equ:gamma}
\left\{\gamma^{\mu},\gamma^{\nu}\right\}:=\gamma^{\mu}\gamma^{\nu}+\gamma^{\nu}\gamma^{\mu}=-2g_{\mu\nu}I_{4}\quad \mbox{and}\quad (\gamma^{\mu})^{*}=-g_{\mu\nu}\gamma^{\nu}.
\end{equation}
Here, $I_{4}$ is the $4\times 4$ identity matrix, $A^{*}=(\bar{A})^{T}$ is the Hermitian conjugate of the matrix $A$ and $g=\rm{diag}\left(-1,1,1,1\right)$ denotes the Minkowski metric in $\R^{1+3}$ with which the indices are raised or lowered. To be more concrete, the Dirac matrices can be represented by
\begin{equation*}
\gamma^{0}=\begin{pmatrix}
I_{2} & 0\\
0 & -I_{2}
\end{pmatrix}\quad \mbox{and}\quad 
\gamma^{a}=\begin{pmatrix}
0 & \sigma_{a}\\
-\sigma_{a} & 0
\end{pmatrix}\quad \mbox{for}\ a=1,2,3,
\end{equation*}
where 
\begin{equation*}
I_{2}=\begin{pmatrix}
1 & 0\\
0 & 1
\end{pmatrix},\quad  
\sigma_{1}=\begin{pmatrix}
0 & 1\\
1 & 0
\end{pmatrix},\quad 
\sigma_{2}=\begin{pmatrix}
0 & -i\\
i & 0
\end{pmatrix},\quad 
\sigma_{3}=\begin{pmatrix}
1 & 0\\
0 & -1
\end{pmatrix}.
\end{equation*}
 We recall that $(V, U) = (I_4, \gamma^0)$ corresponds to the scalar case, while $(V, U) = (i\gamma^5, i\gamma^0 \gamma^5)$ corresponds to the pseudoscalar case where $\gamma^{5}=i\gamma^{0}\gamma^{1}\gamma^{2}\gamma^{3}$.

\smallskip
The main purpose of the present paper is to show the uniform boundedness of the total energy of the Dirac--Klein-Gordon system \eqref{eq:DKG-Mm}. We note that there are eight cases in total if we consider only $m,M \in \{0, 1\}$, and $(V, U) = (I_4, \gamma^0)$ or $(V, U) = (i\gamma^5, i\gamma^0 \gamma^5)$. First, for the case of $(m,M)=(1,1)$, we recall that the system \eqref{eq:DKG-Mm} is reduced to coupled Klein-Gordon equations, and thus the classical results of Klainerman \cite{Klainerman85} and Shatah \cite{Shatah} imply the uniform boundedness of the total energy. Second, for the case of $(m,M)=(1,0)$, thanks to the faster decay rates of the 3D case than the 2D case, we expect that the same strategy used in~\cite{DLMY,DongZoeDKG} can lead to the uniform boundedness for the total energy of~\eqref{eq:DKG-Mm}. Therefore, what is still open is the remaining four cases for $(m,M)=(0,0)$ and $(m,M)=(0,1)$.

\smallskip
$\bullet$ {\textbf{Case I}}: $(m,M,V,U)=(0,0,I_{4},\gamma^{0})$. 
Consider the 3D massless Dirac--Klein-Gordon system 
\begin{equation}\label{equ:3DDKG}
\left\{\begin{aligned}
&-i\gamma^{\mu}\partial_{\mu}\psi=\phi\psi,\quad (t,x)\in [2,\infty)\times \R^{3},\\
&-\Box\phi=\psi^{*}\gamma^{0}\psi,\ \quad (t,x)\in [2,\infty)\times \R^{3}.
\end{aligned}\right.
\end{equation}
The initial data are prescribed at $t=2$
\begin{equation}\label{equ:ID}
\big( \psi, \phi, \partial_t \phi \big)_{|t=2}
=(\psi_{0},\phi_{0},\phi_{1}),
\end{equation}
in which $(\psi_{0},\phi_{0},\phi_{1}): \R^{3}\to \C^{4}\times \R\times \R$.
In particle
physics, the system~\eqref{equ:3DDKG}-\eqref{equ:ID} describes a massless spinor $\psi$ and a massless scalar field $\phi$, whose interaction Lagrangian is given by 
\begin{equation*}
\mathcal{L}_{\rm{Yukawa}}(\psi,\phi)=\psi^{*}\gamma^{0}\phi\psi.
\end{equation*}

\smallskip
$\bullet$  \textbf{Case II}: $(m,M,V,U)=(0,0,i\gamma^{5},i\gamma^{0}\gamma^{5})$.
Consider the following 3D massless Dirac--Klein-Gordon system 
\begin{equation}\label{equ:3DDKGcase2}
\left\{\begin{aligned}
&-i\gamma^{\mu}\partial_{\mu}\psi=i\phi\gamma^{5}\psi,\quad (t,x)\in [2,\infty)\times \R^{3},\\
&-\Box\phi=i\psi^{*}\gamma^{0}\gamma^{5}\psi,\ \quad (t,x)\in [2,\infty)\times \R^{3}.
\end{aligned}\right.
\end{equation}
The initial data are prescribed at $t=2$
\begin{equation}\label{equ:IDcase2}
\big( \psi, \phi, \partial_t \phi \big)_{|t=2}
=(\psi_{0},\phi_{0},\phi_{1}),
\end{equation}
in which $(\psi_{0},\phi_{0},\phi_{1}): \R^{3}\to \C^{4}\times \R\times \R$.
The nonlinear terms of~\eqref{equ:3DDKG} and~\eqref{equ:3DDKGcase2} are different, however they enjoy  some similar properties (see more details in \S\ref{S:thmmassless}). Therefore, we expect that this case can be treated using the same strategy as in Case I.

\smallskip 

$\bullet$ \textbf{Case III}: $(m,M,V,U)=(0,1,i\gamma^{5},i\gamma^{0}\gamma^{5})$.
We also consider the following 3D Dirac--Klein-Gordon system 
\begin{equation}\label{equ:3DDKG2}
\left\{\begin{aligned}
&-i\gamma^{\mu}\partial_{\mu}\psi+\psi=i\phi\gamma^{5}\psi,\quad (t,x)\in [2,\infty)\times \R^{3},\\
&-\Box\phi=i\psi^{*}\gamma^{0}\gamma^{5}\psi, \quad \quad  \quad (t,x)\in [2,\infty)\times \R^{3}.
\end{aligned}\right.
\end{equation}
The initial data are also prescribed at $t=2$
\begin{equation}\label{equ:ID2}
\big( \psi, \phi, \partial_t \phi \big)_{|t=2}
=(\psi_{0},\phi_{0},\phi_{1}),
\end{equation}
in which $(\psi_{0},\phi_{0},\phi_{1}): \R^{3}\to \C^{4}\times \R\times \R$.
In  particle physics, the system~\eqref{equ:3DDKG2}-\eqref{equ:ID2} describes a massive spinor $\psi$ and a massless pseudoscalar field $\phi$, whose interaction Lagrangian is given by 
\begin{equation*}
\mathcal{L}_{\rm{Yukawa}}(\psi,\phi)=i\psi^{*}\gamma^{0}\gamma^{5}\phi\psi.
\end{equation*}

\smallskip
$\bullet$ \textbf{Case IV}: $(m,M,V,U)=(0,1,I_{4},\gamma^{0})$. Consider the following 3D Dirac--Klein-Gordon system 
\begin{equation}\label{equ:3DDKG2case4}
\left\{\begin{aligned}
&-i\gamma^{\mu}\partial_{\mu}\psi+\psi=\phi\psi,\quad (t,x)\in [2,\infty)\times \R^{3},\\
&-\Box\phi=\psi^{*}\gamma^{0}\psi, \quad \quad  \quad (t,x)\in [2,\infty)\times \R^{3}.
\end{aligned}\right.
\end{equation}
The initial data are also prescribed at $t=2$
\begin{equation}\label{equ:ID2case4}
\big( \psi, \phi, \partial_t \phi \big)_{|t=2}
=(\psi_{0},\phi_{0},\phi_{1}),
\end{equation}
in which $(\psi_{0},\phi_{0},\phi_{1}): \R^{3}\to \C^{4}\times \R\times \R$. Note that, the nonlinear term's structures of Case III and Case IV are quite different, and in particular, these structures play an important role in the proof for the Case III (see more details in \S\ref{SS:Keynonlinear}). Therefore, we do not expect that this case can be treated using the same strategy as in Case III, even though the systems of the two cases enjoy some similarities (at least in the linear level). Actually, we are not yet able to treat this last case, and we leave it as an interesting open question.


\subsection{Main results}\label{SS:main}
From now on, we focus on the Case I and Case III listed in the above.
Our first main result is formulated as follows.
\begin{theorem}\label{thm:massless}
	Let $N\ge 7$ be an integer and $\mathcal{X}^{N}=(\dot{H}^{N+1}\cap \dot{H}^{1})\times (\dot{H}^{N}\cap L^{2})$. There exists an $\epsilon_{0}>0$ such that for all $\epsilon\in (0,\epsilon_{0})$, and all compactly supported initial data satisfying the smallness condition
	\begin{equation}\label{est:small}
    \|\psi_{0}\|_{H^{N}}+\|\phi_{0}\|_{H^{N+1}}+\|\phi_{1}\|_{H^{N}}\le \epsilon,
	\end{equation}
	the Cauchy problem~\eqref{equ:3DDKG}-\eqref{equ:ID} admits a global-in-time solution $(\psi,\phi)$ which satisfies the following sharp time decay results
	\begin{equation}\label{est:theorem1point}
	\left|\psi(t,x)\right|\lesssim \epsilon(t+|x|)^{-1}\langle t-|x|\rangle^{-\frac{1}{2}},\quad 
	\left|\phi(t,x)\right|\lesssim \epsilon(t+|x|)^{-1}.
	\end{equation}
	With $\Gamma=\left\{\partial,H,\Omega,S\right\}$, the following uniform energy estimates are also valid
	\begin{equation}\label{est:theorem1}
	\sup_{t\in [2,\infty)}\sum_{|I|\le N}\left(\left\|\Gamma^{I}\psi(t,x)\right\|_{L_{x}^{2}}+\left\|\partial \Gamma^{I}\phi(t,x)\right\|_{L_{x}^{2}}\right)\lesssim \epsilon.
	\end{equation}
	Moreover, the massless Klein-Gordon field $\phi$ scatters to a free solution in $\mathcal{X}^{N}$: there exists $(\phi_{0\ell},\phi_{1\ell})\in \mathcal{X}^{N}$ such that 
	\begin{equation}\label{est:theorem1sca}
	\lim_{t\to \infty}\left\|(\phi,\pt \phi)-(\phi_{\ell},\pt \phi_{\ell})\right\|_{\mathcal{X}^{N}}=0,
	\end{equation}
	where $\phi_{\ell}$ is the solution to the 3D linear homogeneous wave equation with the initial data $(\phi_{0\ell},\phi_{1\ell})$.
\end{theorem}

\begin{remark}
The symbol $\Gamma$ in Theorem \ref{thm:massless} represents vector fields that are compatible with the wave operator, whose definition is given in \S\ref{SS:Notation}.
\end{remark}

Our second main result is formulated as follows.
\begin{theorem}\label{thm:massive}
	Let $N\ge 13$ be an integer and $\mathcal{X}^{N}=(\dot{H}^{N+1}\cap \dot{H}^{1})\times (\dot{H}^{N}\cap L^{2})$. There exists an $\varepsilon_{0}>0$ such that for all $\varepsilon\in (0,\varepsilon_{0})$, and all initial data satisfying the smallness condition
	\begin{equation}\label{est:small2}
	\begin{aligned}
	\sum_{k\le N}\left\|\log (2+|x|)\langle x\rangle^{k+\frac{3}{2}}\nabla_{x}^{k}\psi_{0}\right\|_{L^2_x}&\lesssim \varepsilon,\\
	\sum_{k\le N+1}\left\|\langle x\rangle^{k}\nabla_{x}^{k}\phi_{0}\right\|_{L_{x}^{2}}+\sum_{k\le N}\left\|\langle x\rangle^{k}\nabla_{x}^{k}\phi_{1}\right\|_{L^2_x}&\lesssim \varepsilon,
	\end{aligned}
	\end{equation}
	the Cauchy problem~\eqref{equ:3DDKG2}-\eqref{equ:ID2} admits a global-in-time solution $(\psi,\phi)$ which satisfies the following sharp pointwise decay results
	\begin{equation}\label{est:theorem2point}
	\left|\psi(t,x)\right|\lesssim \varepsilon(t+|x|)^{-\frac{3}{2}},\quad 
	\left|\phi(t,x)\right|\lesssim \varepsilon(t+|x|)^{-1}\langle t-|x|\rangle^{-\frac{1}{2}}.
	\end{equation}
	With $Z=\left\{\partial,H,\Omega\right\}$, the following uniform energy estimates are also valid
	\begin{equation}\label{est:theorem2}
	\sup_{t\in [2,\infty)}\sum_{|I|\le N}\left(\left\|Z^{I}\psi(t,x)\right\|_{L_{x}^{2}}+\left\| Z^{I}\phi(t,x)\right\|_{L_{x}^{2}}+\left\|\partial Z^{I}\phi(t,x)\right\|_{L_{x}^{2}}\right)\lesssim \varepsilon.
	\end{equation}
	Moreover, the solution $(\psi,\phi)$ scatters to a free solution in $H^{N-2}\times\mathcal{X}^{N}$: there exists $(\psi_{0\ell},\phi_{0\ell},\phi_{1\ell})\in H^{N-2}\times \mathcal{X}^{N}$ such that 
	\begin{equation}\label{est:theorem2sca}
	\lim_{t\to \infty}\left\|(\psi,\phi,\pt \phi)-(\psi_{\ell},\phi_{\ell},\pt \phi_{\ell})\right\|_{H^{N-2}\times \mathcal{X}^{N}}=0,
	\end{equation}
	where $\left(\psi_{\ell},\phi_{\ell}\right)$ is the solution to the 3D linear homogeneous Dirac--Klein-Gordon system with the initial data $(\psi_{0\ell},\phi_{0\ell},\phi_{1\ell})$.
\end{theorem}

\begin{remark}
The symbol $Z$ in Theorem \ref{thm:massive} represents vector fields that are compatible with both the wave operator and the Klein-Gordon operator, whose definition can be found in \S\ref{SS:Notation}.
\end{remark}

\begin{remark}
	It is worth mentioning that the uniform boundedness of the element $\left\| Z^{I}\phi(t,x)\right\|_{L_{x}^{2}}$ in \eqref{est:theorem2} is not trivial for waves, as the basic energy for the wave equation can only control the differentiated wave components $\left\| \partial Z^{I}\phi(t,x)\right\|_{L_{x}^{2}}$.
\end{remark}

\begin{remark}
	For $(m,M)=(0,1)$, Bachelot~\cite{Bache} first proved stability and pointwise decay while Theorem \ref{thm:massive} provides a uniform boundedness result for the total energy in the pseudoscalar case, i.e. $(V,U)=(i\gamma^5,i \gamma^0 \gamma^5)$. But it is still open whether such a result is valid in the scalar case, i.e. $(V,U)=(I_{4},\gamma^{0})$. 
\end{remark}

	Recall that, the global existence results to \eqref{equ:3DDKG}-\eqref{equ:ID} and \eqref{equ:3DDKG2}-\eqref{equ:ID2} were first shown by Choquet-Bruhat in \cite{Choquet} using the conformal method of Christodoulou and by Bachelot in~\cite{Bache} using the energy method, respectively. Thus our main contributions to \eqref{equ:3DDKG}-\eqref{equ:ID} and \eqref{equ:3DDKG2}-\eqref{equ:ID2} are the uniform boundeness of the total energies \eqref{est:theorem1}-\eqref{est:theorem2}. The motivations for proving a uniform boundedness result mainly lie in three points: 1) It is natural to explore whether the total energy for a physical model is uniformly bounded or not. If the total energy for the Dirac--Klein-Gordon system can be proved uniformly bounded in time, then such a result is an optimal one in the sense that the linear system has conserved energy. 2) The uniform boundedness result is a necessary result for linear scattering, which is a vital question concerned by the society of dispersive equations. Indeed, in Theorem \ref{thm:massless}, we only provide a linear scattering result for the massless Klein-Gordon field $\phi$, while the one for the Dirac field $\psi$ is absent. 3) A uniform boundedness result in 3D for Dirac--Klein-Gordon system can be regarded as a stepping stone towards the asymptotic stability for the system in 2D, since the decay rates of 2D Dirac, wave and Klein-Gordon equations are slower than those in the 3D case.

\subsection{Brief discussion on related results}

The study on Dirac--Klein-Gordon system has attracted people's attention since decades ago. The pioneering works on its global existence in 3D are due to Chadam and Glassey \cite{ChGl74} for some restricted initial data, and to Choquet-Bruhat \cite{Choquet} in the massless case \eqref{equ:3DDKG}-\eqref{equ:ID}. The first stability result in 3D with pointwise decay was proved by Bachelot \cite{Bache} for \eqref{equ:3DDKG2}-\eqref{equ:ID2}, which is closely related to our study. As for the 2D case, the case $(m,M)=(1,1)$ can be reduced to coupled Klein-Gordon equations, which were tackled by Simon-Taflin \cite{SimTaf}, by Ozawa-Tsutaya-Tsutsumi \cite{OTT}, and by Delort-Fang-Xue \cite{Delort} with pointwise decay result, and for all combinations of masses, Gr\"{u}nrock-Pecher \cite{GH10} first proved the global existence result. Recently, the pointwise decay and scattering for the 2D case with $(m,M)=(1,0)$ were established in \cite{DLMY,DongZoeDKG}. We also remind one of some interesting results regarding low-regular initial data \cite{Bejenaru-Herr,DFS-07, WangX}.

\smallskip
We recall that global existence and pointwise decay for Case IV, i.e. $(m,M,V,U)=(1,0,I_{4},\gamma^{0})$, was demonstrated in \cite{DW-JDE,Ionescu-P-WKG}. But it is still unknown whether the total energy in this case is uniformly bounded in time or not. 

\smallskip
To explore whether the total energy is uniformly bounded in time is of importance to several types of nonlinear evolution equations, one finds such study on incompressible elastodynamics \cite{Cai2020,LeiWan} (see also references therein), for instance. Such studies are closely related to linear scattering, and the motivations are briefly discussed in the end of \S\ref{SS:main} and in \cite{Cai2020,LeiWan}. Last, we mention the results \cite{BMS,Candy-Herr2, Chen-Dirac, DM20, Guo-N-W, Ionescu-P-EKG, Klainerman-WY, PLF-YM-cmp, WangQ} which are also relevant to our study.

\subsection{Key strategies and novelties}

To start, we briefly mention the major difficulties encountered when proving Theorems \ref{thm:massless} and \ref{thm:massive}: nonintegrable decay of wave components in 3D, slow decay of the nonlinearities, lack of scaling invariance when proving Theorem \ref{thm:massive}, lack of derivatives on the massless Klein-Gordon field $\phi$ in the nonlinearities of the Dirac equation. Some of these are not issues if we only want to show global existence and pointwise decay, but they turn into obstacles when we establish a uniform boundedness result. One finds a detailed discussion in recent results \cite{DLMY,DongZoeDKG} for instance. In what follows, we list the three key points in the proof that can be used to overcome the difficulties listed in the above.
 
\smallskip
First, we recall that the Klainerman-Sobolev inequality \eqref{est:Klainer} gives $\langle t+r \rangle^{-1} \langle t-r\rangle^{-1/2}$ decay for the linear wave component. Note that, this rate is a non-integrable quantity in time $t$, which further makes it nontrivial to show a uniform boundedness result for the total energy of \eqref{equ:3DDKG}-\eqref{equ:ID} and \eqref{equ:3DDKG2}-\eqref{equ:ID2}. Fortunately, the ghost weight method by Alinhac \cite{Alinhac01b} (possibly adapted to Klein-Gordon setting) together with some special structure of the Dirac equation and the nonlinearities used in \cite{DLMY,DongZoeDKG} allows one to turn $\langle t-r \rangle$ decay into $\langle t+r \rangle$ decay, which solves the problem of insufficient time decay of the solution $\phi$.

\smallskip
Second, one key issue is how to get sufficient bounds for the undifferentiated wave component $\phi$ (as it appears in the nonlinearities of the Dirac equation), both in the $L^{2}$ and $L^\infty$ levels. The contents of the strategies are split into two parts according to the norms that are studied:

\smallskip
$\bullet$ To show the $L^2$ bounds of $\phi$, in the proof of Theorem \ref{thm:massless}, the key is to establish a $(t-r)$-weighted conformal energy estimate for the wave equation inspired by earlier works \cite{Alinhac01b, Dong2006} (see Lemma \ref{le:weightcon}), which seems to be new. Roughly speaking, this weighted conformal energy estimate allows us to benefit from the $(t-r)/t$ decay in the source term, and in turn the $(t-r)/t$ decay can be inherited by $\phi$. In the proof of Theorem~\ref{thm:massive}, the $L^2$ norm of $\phi$ can be estimated using the conformal energy. 

\smallskip
$\bullet$ To show the $L^\infty$ bounds of $\phi$, in the proof of Theorem \ref{thm:massless}, the analysis is involved. We first apply an observation in \cite{Boura00} (see~\eqref{equ:defPsi}) to transform the Dirac equation of $\psi = -i\gamma^{\mu}\partial_{\mu}\Psi$ into a wave equation of $\Psi$. We then apply an $L^\infty-L^\infty$ estimate of \cite{AlinhacIndiana} (see Proposition \ref{pro:Linfinity}) to carefully derive the pointwise decay of $\phi$.
In the proof of Theorem \ref{thm:massive} and thanks to a hidden nonlinear transformation of \cite{Bache} (see Lemma \ref{le:Hidden}), we can direct apply an $L^\infty-L^\infty$ estimate of \cite{Kubota-Yokoyama} (see Proposition \ref{pro:Linwave}) to get desired pointwise decay of $\phi$.

\smallskip
Third, we cannot commute with the scaling vector field in the proof of Theorem \ref{thm:massive} due to the presence of the mass term in the Dirac equation. In this case, a weaker version of global Sobolev inequality \eqref{est:SobolevGlobal} is used to obtain rough decay results for the solution, and then the decay results can be improved by an $L^\infty-L^\infty$ estimate of \cite{Kubota-Yokoyama} (see Proposition \ref{pro:Linwave}) and a pointwise estimate of \cite{Geor} (see Proposition \ref{prop:decayKG}).

\subsection{Organization} The article is organized as follows. In Section~\ref{S:pre}, we introduce
the notations and some fundamental analysis tools to be used in Sections~\ref{S:thmmassless}-\ref{S:thmmassive}: the estimates on commutators and null forms, global Sobolev inequality, and energy and $L^{\infty}$ estimates for the Dirac, wave and Klein-Gordon equations. In Sections~\ref{S:thmmassless}-\ref{S:thmmassive}, by Klainerman’s vector field method enhanced with
Alinhac’s ghost weight method, we provide the proof for Theorem~\ref{thm:massless} and Theorem~\ref{thm:massive}, respectively.

\section*{Acknowledgment}
The author S. D. would like to thank Zoe Wyatt (KCL) for helpful discussions.
 Part of this work was done while the author X.Y. was visiting the Institute of Mathematics, AMSS, Chinese Academy of Sciences, and X.Y. thanks the institution for its warm hospitality.

\section{Preliminaries}\label{S:pre}
\subsection{Notations and conventions}\label{SS:Notation}
We work in the (1+3) dimensional Minkowski spacetime $\R^{1+3}$. We denote one point by $(t,x)=(x_{0},x_{1},x_{2},x_{3})$ and its spatial radius by $r=\sqrt{x_{1}^{2}+x_{2}^{2}+x_{3}^{2}}$. Denote $\omega_{a}=\frac{x_{a}}{r}$ for $a=1,2,3$ and $x=(x_{1},x_{2},x_{3})\in \R^{3}$.
Greek indices $\left\{\mu,\nu,\dots\right\}$ range over $\left\{0,1,2,3\right\}$, Roman indices $\left\{a,b,\dots\right\}$ range over $\left\{1,2,3\right\}$, and the Einstein summation convention for repeated upper and lower indices is applied throughout this paper.
We use $\Box = \partial_\alpha \partial^\alpha = -\partial_0^{2}+ \partial_{1}^{2}+\partial_{2}^{2}+\partial_{3}^{2}$ to denote the d'Alembert operator. 

To state global
Sobolev inequalities, we first introduce the following four groups of vector fields:
\begin{enumerate}
	
		\smallskip
	\item [(i)]  Scaling vector field: $S=t\partial_{t}+x^{a}\partial_{a}$.
	
	\smallskip
	\item [(ii)] Translations: $\partial_{\alpha}=\partial_{x_{\alpha}}$, for $\alpha=0,1,2,3$.
	
	\smallskip
	\item [(iii)] Hyperbolic rotations: $H_{a}=x_{a}\partial_{t}+t\partial_{a}$, for $a=1,2,3$.
	
	\smallskip
	\item [(iv)] Spatial rotations: $\Omega_{ab}=x_{a}\partial_{b}-x_{b}\partial_{a}$, for $1\le a<b\le 3$.
\end{enumerate}

We denote 
\begin{equation*}
\partial=\left(\partial_{0},\partial_{1},\partial_{2},\partial_{3}\right),\ 
H=\left(H_{1},H_{2},H_{3}\right)\ \mbox{and}\
\Omega=\left(\Omega_{12},\Omega_{13},\Omega_{23}\right).
\end{equation*}

Following Bachelot~\cite{Bache}, we also introduce the following modified spatial rotations and hyperbolic rotations which commute with the Dirac operator,
\begin{equation*}
\begin{aligned}
\widehat{H}_{a}&=H_{a}-\frac{1}{2}\gamma^{0}\gamma^{a},\quad \widehat{H}=\big(\widehat{H}_{1},\widehat{H}_{2},\widehat{H}_{3}\big),\\
\widehat{\Omega}_{ab}&=\Omega_{ab}-\frac{1}{2}\gamma^{a}\gamma^{b},\quad \widehat{\Omega}=\big(\widehat{\Omega}_{12},\widehat{\Omega}_{13},\widehat{\Omega}_{23}\big).
\end{aligned}
\end{equation*}
We define the following four ordered sets of vector fields,
\begin{equation*}
\begin{aligned}
Z:=\left(Z_{1},\dots,Z_{10}\right)&=\left(\partial,H,\Omega\right),\quad \Gamma:=\left(\Gamma_{1},\dots,\Gamma_{11}\right)=\left(S,\partial,H,\Omega\right),\\
\widehat{Z}:=\big(\widehat{Z}_{1},\dots,\widehat{Z}_{10}\big)&=\big(\partial,\widehat{H},\widehat{\Omega}\big),\quad 
\widehat{\Gamma}:=\big(\widehat{\Gamma}_{1},\dots,\widehat{\Gamma}_{11}\big)=\big(S,\partial,\widehat{H},\widehat{\Omega}\big).
\end{aligned}
\end{equation*}
Moreover, for any $\R$-valued or $\C^{4}$-valued regular function $f=f(t,x)$, we denote 
\begin{equation*}
\begin{aligned}
\left|Zf\right|&=\left(\sum_{k=1}^{10}|Z_{k}f|^{2}\right)^{\frac{1}{2}},\quad 
\left|\Gamma f\right|=\left(\sum_{k=1}^{11}|\Gamma_{k}f|^{2}\right)^{\frac{1}{2}},\\
\big|\widehat{Z}f\big|&=\left(\sum_{k=1}^{10}\big|\widehat{Z}_{k}f\big|^{2}\right)^{\frac{1}{2}},\quad 
\big|\widehat{\Gamma} f\big|=\left(\sum_{k=1}^{11}\big|\widehat{\Gamma}_{k}f\big|^{2}\right)^{\frac{1}{2}}.
\end{aligned}
\end{equation*}
For all multi-index $I=(i_{1},\dots,i_{10})\in \mathbb{N}^{10}$ and $J=(j_{1},\dots,j_{11})\in \mathbb{N}^{11}$, we set 
\begin{align*}
Z^{I}=\prod_{k=1}^{10}Z_{k}^{i_{k}},\quad \widehat{Z}^{I}=\prod_{k=1}^{10}\widehat{Z}_{k}^{i_{k}},\quad 
\Gamma^{J}=\prod_{k=1}^{11}\Gamma_{k}^{j_{k}},\quad \widehat{\Gamma}^{J}=\prod_{k=1}^{11}\widehat{\Gamma}_{k}^{j_{k}}.
\end{align*}
Moreover, we introduce the good derivatives $G_{a}=\partial_{a}+\omega_{a}\partial_{0}$ for $a=1,2,3$ appearing in ghost weight energy estimate.

For all $\R$-valued or $\C^{4}$-valued regular functions $f=f(t,x)$ and $g=g(t,x)$, we consider the standard null forms $Q_{0}$ and $Q_{\alpha\beta}$,
\begin{equation*}
Q_{0}(f,g)=\left(\partial_{\alpha}f\right)^{*}\partial^{\alpha}g\quad \mbox{and}\quad 
Q_{\alpha\beta}(f,g)=\left(\partial_{\alpha}f\right)^{*}\partial_{\beta}g-\left(\partial_{\beta}f\right)^{*}\partial_{\alpha}g.
\end{equation*}

Following~\cite[Section 3]{DongZoeDKG}, we consider the following two terms which are related to the hidden structure of the Dirac$\times$Dirac interaction terms in~\eqref{equ:3DDKG}-\eqref{equ:3DDKG2}: for all regular function $\varphi=\varphi(t,x):\R^{1+3}\to \C^{4}$, we denote 
\begin{equation}\label{equ:def+-}
\begin{aligned}
\left[\varphi\right]_{+}(t,x):&=\left(I_{4}+\omega_{a}\gamma^{0}\gamma^{a}\right)\varphi(t,x),\\
\left[\varphi\right]_{-}(t,x):&=\left(I_{4}-\omega_{a}\gamma^{0}\gamma^{a}\right)\varphi(t,x).
\end{aligned}
\end{equation}

Let $\left\{\zeta_{k}\right\}_{k=0}^{\infty}$ be a Littlewood-Paley partition of unity, i.e.
\begin{equation*}
1=\sum_{k=0}^{\infty}\zeta_{k}(\tau),\ \tau\ge 0,\ \zeta_{k}\in C_{0}^{\infty}\left(\R\right),\quad \zeta_{k}\ge 0\quad \mbox{for all}\ k\ge 0,
\end{equation*}
as well as 
\begin{equation*}
\mbox{supp}\,\zeta_{0}\cap [0,\infty)=[0,2],\quad 
\mbox{supp}\,\zeta_{k}\subset \left[2^{k-1},2^{k+1}\right]\quad \mbox{for all}\ k\ge 1.
\end{equation*}

For simplicity of notation, we denote the initial data $(\psi,\phi_{0},\phi_{1})$ by $(\psi_{0},\vec{\phi}_{0})$. 

Without loss of generality, in the proof of Theorem~\ref{thm:massless} and related technical Lemma of compactly supported functions, we restrict our study to functions supported within the spacetime region $\mathcal{K}:=\left\{(t,x)\in [2,\infty)\times \R^{3}:r+1\le t\right\}$.

\subsection{Basic estimates}

In this subsection, we recall some preliminary estimates related to commutators and vector fields.
First, since the differences between vector fields $(Z,\Gamma)$ and modified vector fields
$\big(\widehat{Z},\widehat{\Gamma}\big)$ are constant matrices, we obtain
\begin{equation}\label{est:GammaGammahat}
\begin{aligned}
\sum_{|I|\le K}\big|\Gamma^{I}f\big|\lesssim \sum_{|I|\le K}\big|\widehat{\Gamma}^{I}f\big|\lesssim \sum_{|I|\le K}\big|\Gamma^{I}f\big|,\\
\sum_{|I|\le K}\big|Z^{I}f\big|\lesssim \sum_{|I|\le K}\big|\widehat{Z}^{I}f\big|\lesssim \sum_{|I|\le K}\big|Z^{I}f\big|,
\end{aligned}
\end{equation}
\begin{equation}\label{est:pGammaGammahat}
\begin{aligned}
\sum_{|I|\le K}\big|\partial\Gamma^{I}f\big|\lesssim \sum_{|I|\le K}\big|\partial\widehat{\Gamma}^{I}f\big|\lesssim \sum_{|I|\le K}\big|\partial\Gamma^{I}f\big|,\\
\sum_{|I|\le K}\big|\partial Z^{I}f\big|\lesssim \sum_{|I|\le K}\big|\partial \widehat{Z}^{I}f\big|\lesssim \sum_{|I|\le K}\big|\partial Z^{I}f\big|,
\end{aligned}
\end{equation}
for any regular $\R$-valued or $\mathbb{C}^{4}$-valued function $f$ and ${K}\in \mathbb{N}^{+}$. Note that, from the definition of $\big(Z,\Gamma,\widehat{Z},\widehat{\Gamma}\big)$, we have the following identities:
\begin{equation}\label{equ:comm}
\begin{aligned}
\left[-\Box,Z_{k}\right]=0\quad &\mbox{and}\quad \left[-\Box,S\right]=-2\Box,\\
\big[-i\gamma^{\mu}\partial_{\mu},\widehat{Z}_{k}\big]=0\quad &\mbox{and}\quad  [-i\gamma^{\mu}\partial_{\mu},S]=-i\gamma^{\mu}\partial_{\mu}.
\end{aligned}
\end{equation}

Second, we recall the following estimates related to the vector fields from~\cite{Sogge}.
\begin{lemma}[\cite{Sogge}]
	For all regular function $f=f(t,x)$, the following hold.
	\begin{enumerate}
		\item \emph{Estimates on commutators}. For all $I\in \mathbb{N}^{10}$ and $J\in \mathbb{N}^{11}$, we have 
		\begin{equation}\label{est:parGamma}
		\begin{aligned}
		\sum_{k=1}^{10}\left|\left[Z_{k},Z^{I}\right]f\right|&\lesssim \sum_{|K|\le |I|}\left|Z^{K}f\right|,\\
		\left|\left[\partial,Z^{I}\right]f\right|
		+\left|\left[S,Z^{I}\right]f\right|&\lesssim \sum_{|K|<|I|}\left|\partial Z^{K}f\right|,\\		\left|\left[\partial,\Gamma^{J}\right]f\right|
		+\left|\left[S,\Gamma^{J}\right]f\right|&\lesssim \sum_{|K|<|J|}\left|\partial \Gamma^{K}f\right|.
		\end{aligned}
		\end{equation}
		
		\item \emph{Estimates on $\partial f$ and $G_{a}f$.} We have 
		\begin{equation}\label{est:Ga}
		\langle t-r\rangle\left|\partial f\right|+\langle t+r\rangle\sum_{a=1}^{3}\left|G_{a}f\right|\lesssim \left|\Gamma f\right|.
		\end{equation}
	\end{enumerate}
\end{lemma}

\begin{proof}
	We refer to~\cite[Page 38 and Proposition 1.1]{Sogge} for the details of the proof.
	\end{proof}

Third, we recall the following estimates related to the null form $Q_{\alpha\beta}$ from~\cite{Sogge}.

\begin{lemma}[\cite{Sogge}]
		Let $f=f(t,x)$ and $g=g(t,x)$ be two regular functions on $\R^{1+3}$, then for any multi-index $I\in \mathbb{N}^{10}$, we have 
		\begin{equation}\label{est:Qab}
		\sum_{\alpha \ne \beta}\left|Z^{I}Q_{\alpha\beta}(f,g)\right|\lesssim \langle t+r\rangle^{-1}\sum_{|I_{1}|+|I_{2}|\le I}
		\left|ZZ^{I_{1}}f\right|\left|ZZ^{I_{2}}g\right|.
		\end{equation}
\end{lemma}

\begin{proof}
	We refer to~\cite[Lemma 3.3]{Sogge} and~\cite[Page 58]{Sogge} for the details of the proof.
	\end{proof}

Next, we recall the following local and global Sobolev inequalities.
\begin{theorem}[\cite{Geor,KlainWave}]
	For all regular function $f=f(t,x)$, the following hold.
	\begin{enumerate}
		\item \emph{Klainerman-Sobolev inequality.} We have 
	\begin{equation}\label{est:Klainer}
	|f(t,x)|\lesssim \langle t+r\rangle^{-1}\langle t-r\rangle^{-\frac{1}{2}}\sum_{|I|\le 2}\left\|\Gamma^{I}f(t,x)\right\|_{L^{2}_{x}}.
	\end{equation}
	
	\item \emph{Standard Sobolev inequality}. We have 
	\begin{equation}\label{est:SobolevStand}
	\left|f(t,x)\right|\lesssim \langle r \rangle^{-1}\sum_{|I|\le 2}\left\|Z^{I}f(t,x)\right\|_{L^{2}_{x}}.
	\end{equation}
	
	\item \emph{Estimate inside of a cone.} For $|x|\le \frac{t}{2}$, we have 
	\begin{equation}\label{est:Sobolevinsi}
	\left|f(t,x)\right|\lesssim \langle t\rangle^{-\frac{3}{4}}\sum_{|I|\le 3}\left\|Z^{I}f(t,x)\right\|_{L^{2}_{x}}.
	\end{equation}
	
	\item \emph{Global Sobolev inequality.} We have 
	\begin{equation}\label{est:SobolevGlobal}
	\left|f(t,x)\right|\lesssim \langle t+r\rangle^{-\frac{3}{4}}\sum_{|I|\le 3}\left\|Z^{I}f(t,x)\right\|_{L_{x}^{2}}.
	\end{equation}
	
\end{enumerate}
\end{theorem}

\begin{proof}
		For~\eqref{est:Klainer}-\eqref{est:Sobolevinsi}, we refer to~\cite[Proposition 1-2]{KlainWave} and~\cite[Lemma 2.4]{Geor} for the details of the proof respectively. Then, the inequality~\eqref{est:SobolevGlobal} is a direct consequence of~\eqref{est:SobolevStand} and~\eqref{est:Sobolevinsi}.
	\end{proof}

Then, we recall the following refined Hardy inequality from~\cite{LindCMAP}.
\begin{lemma}[\cite{LindCMAP}]
	For all regular function $f=f(t,x)$ with support in $\mathcal{K}$, we have 
	\begin{equation}\label{est:Hardy}
	\left\|\frac{f(t,x)}{(t-r)}\right\|_{L_{x}^{2}}\lesssim \left\|\partial_{r} f(t,x)\right\|_{L_{x}^{2}}.
	\end{equation}
\end{lemma}
\begin{proof}
	We refer to~\cite[Lemma 1.2]{LindCMAP} for the details of the proof.
	\end{proof}

Last, we introduce the following technical estimates related to the Dirac$\times$Dirac interaction term which will appear in the energy estimate.

\begin{lemma}
		Let $\Theta_{1}$ be an $\R$-valued regular function and $\Theta_{2}$ be a $\C^{4}$-valued regular function on $\R^{1+3}$, then for any multi-index $I\in \mathbb{N}^{10}$ and $J\in \mathbb{N}^{11}$, we have 
		\begin{equation}\label{est:Z-}
		\big|\big[\widehat{Z}^{I}\big(\Theta_{1}\Theta_{2}\big)\big]_{-}\big|\lesssim \sum_{|I_{1}|+|I_{2}|\le I}\big|Z^{I_{1}}\Theta_{1}\big|\big|\big[\widehat{Z}^{I_{2}}\Theta_{2}\big]_{-}\big|,
		\end{equation}
		\begin{equation}\label{est:Gamma-}
		\big|\big[\widehat{\Gamma}^{J}\big(\Theta_{1}\Theta_{2}\big)\big]_{-}\big|\lesssim \sum_{|J_{1}|+|J_{2}|\le J}\big|\Gamma^{J_{1}}\Theta_{1}\big|\big|\big[\widehat{\Gamma}^{J_{2}}\Theta_{2}\big]_{-}\big|.
		\end{equation}
\end{lemma}

\begin{proof}
	By an elementary computation, for all $k\in\left\{1,\dots,10\right\}$, we have 
	\begin{equation*}
	\begin{aligned}
	\left[S\left(\Theta_{1}\Theta_{2}\right)\right]_{-}&=\left(S\Theta_{1}\right)\left[\Theta_{2}\right]_{-}
	+\Theta_{1}\left[S\Theta_{2}\right]_{-},\\
	\big[\widehat{Z}_{k}\left(\Theta_{1}\Theta_{2}\right)\big]_{-}&=\left(Z_{k}\Theta_{1}\right)\left[\Theta_{2}\right]_{-}
	+\Theta_{1}\big[\widehat{Z}_{k}\Theta_{2}\big]_{-}.
	\end{aligned}
	\end{equation*}
	Combining the above identities with the Leibniz rule of $\partial$ and then using an induction argument, we obtain~\eqref{est:Z-} and~\eqref{est:Gamma-}.
	\end{proof}

\begin{lemma}
	Let $\Theta_{1}$ and $\Theta_{2}$ be two $\C^{4}$-valued regular functions on $\R^{1+3}$, then the following estimates for $\Theta_{1}^{*}\gamma^{0}\Theta_{2}$ are true.
	\begin{enumerate}
		
	\item \emph{First estimate.} We have 
	\begin{equation}\label{est:Hidden}
	\left|\Theta_{1}^{*}\gamma^{0}\Theta_{2}\right|\lesssim \left|[\Theta_{1}]_{-}\right|\left|\Theta_{2}\right|
	+\left|\Theta_{1}\right|\left|[\Theta_{2}]_{-}\right|.
	\end{equation}
	
	\smallskip
		\item \emph{Second estimate.} For any multi-index $I\in \mathbb{N}^{10}$, we have 
	\begin{equation}\label{est:Z0}
	\left|Z^{I}\left(\Theta_{1}^{*}\gamma^{0}\Theta_{2}\right)\right|\lesssim \sum_{|I_{1}|+|I_{2}|\le |I|}
	\big|\big(\widehat{Z}^{I_{1}}\Theta_{1}\big)^{*}\gamma^{0}\big(\widehat{Z}^{I_{2}}\Theta_{2}\big)\big|.
	\end{equation}

	\item \emph{Third estimate.} For any multi-index $J\in \mathbb{N}^{11}$, we have 
	\begin{equation}\label{est:Gaga0}
\left|\Gamma^{J}\left(\Theta_{1}^{*}\gamma^{0}\Theta_{2}\right)\right|\lesssim \sum_{|J_{1}|+|J_{2}|\le |J|}
\big|\big(\widehat{\Gamma}^{J_{1}}\Theta_{1}\big)^{*}\gamma^{0}\big(\widehat{\Gamma}^{J_{2}}\Theta_{2}\big)\big|.
\end{equation}
	
\end{enumerate}
\end{lemma}

\begin{proof}
	Proof of (i). From the definition of $[\cdot]_{+}$ and $[\cdot]_{-}$ in~\eqref{equ:def+-}, we have 
	\begin{equation*}
	2\Theta_{1}=[\Theta_{1}]_{+}+[\Theta_{1}]_{-}\quad \mbox{and}\quad 
	2\Theta_{2}=[\Theta_{2}]_{+}+[\Theta_{2}]_{-},
	\end{equation*}
	which implies
	\begin{equation*}
	\Theta_{1}^{*}\gamma^{0}\Theta_{2}=\frac{1}{4}\left([\Theta_{1}]_{-}^{*}\gamma^{0}[\Theta_{2}]_{-}+[\Theta_{1}]_{-}^{*}\gamma^{0}[\Theta_{2}]_{+}+[\Theta_{1}]_{+}^{*}\gamma^{0}[\Theta_{2}]_{-}+[\Theta_{1}]_{+}^{*}\gamma^{0}[\Theta_{2}]_{+}\right).
	\end{equation*}
	Note that, for the last term in the above identity, from~\eqref{equ:gamma}, we have 
	\begin{equation*}
	\begin{aligned}
	[\Theta_{1}]_{+}^{*}\gamma^{0}[\Theta_{2}]_{+}
	=&\left(\Theta_{1}^{*}+\omega_{a}\Theta_{1}^{*} \gamma^{0}\gamma^{a}\right)\gamma^{0}\left(\Theta_{2}+\omega_{b} \gamma^{0}\gamma^{b}\Theta_{2}\right)\\
	=&\left(1-\omega_{1}^{2}-\omega_{2}^{2}-\omega_{3}^{2}\right)\Theta_{1}^{*}\gamma^{0}\Theta_{2}+\omega_{a}\Theta_{1}^{*}\gamma^{0}\left(\gamma^{0}\gamma^{a}+\gamma^{a}\gamma^{0}\right)\Theta_{2}\\
	&+\sum_{a<b}\omega_{a}\omega_{b}\Theta_{1}^{*}\gamma^{0}\left(\gamma^{a}\gamma^{b}+\gamma^{b}\gamma^{a}\right)\Theta_{2}=0.
	\end{aligned}
	\end{equation*}
Combining the above identities with the triangle inequality, we obtain~\eqref{est:Hidden}.
	
	Proof of (ii) and (iii). From~\eqref{equ:gamma} and the definition of $\widehat{\Omega}_{ab}$ and $\widehat{H}_{a}$, we have
	\begin{equation*}
	\Omega_{ab}\left(\Theta_{1}^{*}\right)=\big(\widehat{\Omega}_{ab}\Theta_{1}\big)^{*}-\frac{1}{2}\Theta_{1}^{*}\gamma^{a}\gamma^{b}\quad \mbox{and}\quad 
	H_{a}\left(\Theta_{1}^{*}\right)=\big(\widehat{H}_{a}\Theta_{1}\big)^{*}+\frac{1}{2}\Theta_{1}^{*}\gamma^{0}\gamma^{a},
	\end{equation*}
	which implies
		\begin{equation*}
	\begin{aligned}
	\Omega_{ab}\left(\Theta_{1}^{*}\gamma^{0}\Theta_{2}\right)
	&=\big(\widehat{\Omega}_{ab}\Theta_{1}\big)^{*}\gamma^{0}\Theta_{2}+\Theta_{1}^{*}\gamma^{0}\big({\widehat{\Omega}_{ab}}\Theta_{2}\big)-\frac{1}{2}\Theta_{1}^{*}\left(\gamma^{a}\gamma^{b}\gamma^{0}-\gamma^{0}\gamma^{a}\gamma^{b}\right)\Theta_{2}\\
	&=\big(\widehat{\Omega}_{ab}\Theta_{1}\big)^{*}\gamma^{0}\Theta_{2}+\Theta_{1}^{*}\gamma^{0}\big({\widehat{\Omega}_{ab}}\Theta_{2}\big),
	\end{aligned}
	\end{equation*}
	\begin{equation*}
	\begin{aligned}
	H_{a}\left(\Theta_{1}^{*}\gamma^{0}\Theta_{2}\right)
	&=\big(\widehat{H}_{a}\Theta_{1}\big)^{*}\gamma^{0}\Theta_{2}+\Theta_{1}^{*}\gamma^{0}\big({\widehat{H}_{a}}\Theta_{2}\big)+\frac{1}{2}\Theta_{1}^{*}\left(\gamma^{0}\gamma^{a}\gamma^{0}+\gamma^{0}\gamma^{0}\gamma^{a}\right)\Theta_{2}\\
	&=\big(\widehat{H}_{a}\Theta_{1}\big)^{*}\gamma^{0}\Theta_{2}+\Theta_{1}^{*}\gamma^{0}\big({\widehat{H}_{a}}\Theta_{2}\big).
	\end{aligned}
	\end{equation*}
	Combining the above identities with the Leibniz rule of $\partial$ and then using an induction argument, we obtain~\eqref{est:Z0} and~\eqref{est:Gaga0}.
\end{proof}

\subsection{Estimate for the 3D linear Dirac spinor}
In this subsection, we review the ghost weight energy estimate for the Dirac spinor and extra pointwise estimate for the radial derivative of the Dirac spinor. Following~\cite{Alinhac01b, AlinhacBook}, for all $\varphi=\varphi(t,x):\R^{1+3}\to \mathbb{C}^{4}$ and $0<\delta\ll 1$, we consider the ghost weight energy for $\varphi$:
\begin{equation*}
\mathcal{E}^{\delta}_{D}(t,\varphi)=\int_{\R^{3}}\left|\varphi(t,x)\right|^{2}\d x 
+\int_{2}^{t}\int_{\R^{3}}\frac{\left|[\varphi]_{-}(\tau,x)\right|^{2}}{\langle \tau-r\rangle^{1+2\delta}}\d x \d \tau.
\end{equation*}
\begin{lemma}\label{le:DiracGhotst}
	Let $m=0$ or $1$ and $\varphi$ be the solution to the Cauchy problem
	\begin{equation}\label{equ:Diracenergy}
	-i\gamma^{\mu}\partial_{\mu}\varphi+m\varphi=G\quad \mbox{with}\quad \varphi_{|t=2}=\varphi_{0}.
	\end{equation}
	Then we have 
	\begin{equation*}
	\mathcal{E}^{\delta}_{D}(t,\varphi)\lesssim \mathcal{E}^{\delta}_{D}(2,\varphi) +\int_{2}^{t}\int_{\R^{3}}\left|\varphi^{*}(\tau,x)\gamma^{0}G(\tau,x)\right|\d x\d \tau.
	\end{equation*}
\end{lemma}

\begin{proof}
	For $(t,x)\in [2,\infty)\times \R^{3}$, we consider
	\begin{equation*}
	q(t,x)=\int_{-\infty}^{r-t}\langle \tau\rangle^{-1-2\delta}\d \tau,\quad \mbox{where}\ r=|x|.
	\end{equation*}
	Note that, multiplying $i e^q  \varphi^* \gamma^0$ to both sides of the equation~\eqref{equ:Diracenergy}, we have
	\begin{equation} \label{eq:ghostdirac1}
	e^q \varphi^* \partial_0 \varphi + e^q \varphi^*\gamma^0\gamma^a \partial_a \varphi+ime^{q}\varphi^{*}\gamma^{0}\varphi= i e^q \varphi^*\gamma^0G.
	\end{equation}
	Taking the complex conjugate of \eqref{eq:ghostdirac1} and then using~\eqref{equ:gamma}, we also have
	\begin{equation} \label{eq:ghostdirac2}
	e^q (\partial_0 \varphi)^* \varphi + e^q (\partial_a \varphi)^*\gamma^0\gamma^a \varphi-ime^{q}\varphi^{*}\gamma^{0}\varphi = -i e^q G^*\gamma^0 \phi.
	\end{equation}
	Gathering \eqref{eq:ghostdirac1} and \eqref{eq:ghostdirac2}, one can obtain
	\begin{align*}
	\partial_0\left(e^q|\varphi|^2\right) + \partial_a(e^q\varphi^*\gamma^0\gamma^a \varphi) +  \frac{e^q}{\langle r-t\rangle^{1+2\delta}} (|\varphi|^{2} - \omega_a \varphi^*\gamma^0\gamma^a\varphi) =-2e^q {\rm{Im}}(\varphi^*\gamma^0G).
	\end{align*}
	Note also that, from~\eqref{equ:gamma} and the definition of $[\cdot]_{-}$ in~\eqref{equ:def+-}, 
	\begin{align*}
	|[\varphi]_{-}|^2 & = (\varphi-\omega_a \gamma^0\gamma^a \varphi)^*(\varphi-\omega_b\gamma^0\gamma^b \varphi) \\
	&=\left(\varphi^{*}-\omega_{a}\varphi^{*}\gamma^{0}\gamma^{a}\right)(\varphi-\omega_b\gamma^0\gamma^b \varphi)\\
	& =  |\varphi|^{2} - 2\omega_a \varphi^* \gamma^0\gamma^a \varphi + \omega_a\omega_b \varphi^*\gamma^{0}\gamma^{a}\gamma^{0}\gamma^{b} \varphi.
	\end{align*}
	By an elementary computation and~\eqref{equ:gamma}, we compute
	\begin{equation*}
	\omega_a\omega_b \varphi^*\gamma^{0}\gamma^{a}\gamma^{0}\gamma^{b} \varphi=|\varphi|^{2}-\sum_{a<b}\varphi^{*}\left(\gamma^{a}\gamma^{b}+\gamma^{b}\gamma^{a}\right)\varphi=|\varphi|^{2}.
	\end{equation*}
	Gathering the above identities, one can obtain
	\begin{align*} 
	\partial_0(e^q|\varphi|^2) + \partial_a(e^q\varphi^*\gamma^0\gamma^a \varphi) +  \frac{e^q}{2\langle r-t\rangle^{1+2\delta}} \left|[\varphi]_{-}\right|^2 =-2e^q {\rm{Im}}(\varphi^*\gamma^0G).
	\end{align*}
	Last, integrating the above identity over the domain $[2,t]\times \R^{3}$ and then using the fact that $e^{q}\sim 1$, we complete the proof of Lemma~\ref{le:DiracGhotst}.
	\end{proof}

Second, we introduce the extra decay for the radial derivative of the Dirac spinor.
\begin{lemma}\label{le:ptvarphi}
	Let $\varphi$ be the solution to the Cauchy problem
	\begin{equation}\label{equ:Dirac}
	-i\gamma^{\mu}\partial_{\mu}\varphi=G\quad \mbox{with}\quad \varphi_{|t=2}=\varphi_{0},
	\end{equation}
	and assume that $\varphi$ is supported in $\mathcal{K}$. Then we have 
	\begin{equation}\label{est:ptvarphi-}
	\left|\partial_{r}\left( [\varphi]_{-}\right)\right|\lesssim \frac{1}{t+r}
	\big(\big|\widehat{\Gamma} \varphi\big|+|\varphi|\big)+|G|.
	\end{equation}
\end{lemma}

\begin{remark}
 This result provides one with some extra decay of $(t-r)/t$ for the component $ [\varphi]_{-}$, which is favorable when close to the light cone $\{ t=r+1 \}$. To achieve this, one integrates $\partial_{r}\left( [\varphi]_{-}\right)$ along the line connected by the points $(t, t-1)$ and $(t, r_{0})$ for any $0\leq r_0\leq t-1$ (see more details in Lemma~\ref{le:Firstpsi}).
\end{remark}

\begin{proof}[Proof of Lemma~\ref{le:ptvarphi}]
	\textbf{Step 1.} Preliminary estimate. We claim that 
		\begin{equation}\label{est:ptvarphi1}
	\left|\partial_{t}\varphi\right|\lesssim \frac{1}{t-r}\left|\Gamma \varphi\right|+\frac{t+r}{t-r}|G|.
	\end{equation}
	Indeed, from $\partial_{a}=t^{-1}H_{a}-(x_{a}/t)\pt $, we rewrite~\eqref{equ:Dirac} as 
	\begin{equation*}
	\left(\gamma^{0}-\frac{x_{a}}{t}\gamma^{a}\right)\partial_{t}\varphi=-\frac{\gamma^{a}}{t}H_{a}\varphi+iG.
	\end{equation*}
	Multiplying both sides of the above identity by $\left(\gamma^{0}-\left(x_{a}/t\right)\gamma^{a}\right)$ and then using~\eqref{equ:gamma}, we compute
	\begin{equation*}
	\frac{(t+r)(t-r)}{t^{2}}\partial_{t}\varphi=-\frac{1}{t}\left(\gamma^{0}-\frac{x_{a}}{t}\gamma^{a}\right)\gamma^{b}H_{b}\varphi+i\left(\gamma^{0}-\frac{x_{a}}{t}\gamma^{a}\right)G.
	\end{equation*}
	Combining the above identity with ${\rm{supp}}\,\varphi\subseteq \mathcal{K}$, we obtain~\eqref{est:ptvarphi1}.
	
	\smallskip
	\textbf{Step 2.} Conclusion. Using again $\partial_{a}=t^{-1}H_{a}-(x_{a}/t)\pt $, we rewrite~\eqref{equ:Dirac} as
	\begin{equation*}
	\left(\gamma^{0}-\frac{x^{a}}{r}\gamma^{a}\right)\partial_{t}\varphi
	=\gamma^{a}\left(\frac{x_{a}}{t}-\frac{x_{a}}{r}\right)\pt \varphi-\frac{\gamma^{a}}{t}H_{a}\varphi+iG,
	\end{equation*}
	which implies
	\begin{equation*}
\gamma^{0}\partial_{t}\left( [\varphi]_{-}\right)=\frac{\gamma^{a}}{t}\frac{x_{a}}{r}\left[(t-r)\pt \varphi\right]-\frac{\gamma^{a}}{t}H_{a}\varphi+iG.
	\end{equation*}
	It follows from~\eqref{est:ptvarphi1} that
	\begin{equation}\label{est:ptvarphi}
	\left|\partial_{t}\left( [\varphi]_{-}\right)\right|
	\lesssim \frac{1}{t+r}\left|\Gamma \varphi\right|+|G|.
	\end{equation}
	Next, using again $\partial_{a}=t^{-1}H_{a}-(x_{a}/t)\pt $ and $\partial_{r}=(x^{a}/r)\partial_{a}$, we find
	\begin{equation*}
	\frac{x^{a}}{r}H_{a}=t\partial_{r}+r\partial_{t}\Rightarrow
	\partial_{r}=\frac{x^{a}}{r}\frac{H_{a}}{t}-\frac{r}{t}\partial_{t}.
	\end{equation*}
	Therefore, from~\eqref{est:ptvarphi}, it holds
	\begin{equation}\label{est:pavarphi}
	\begin{aligned}
	\left|\partial_{r}\left( [\varphi]_{-}\right)\right|
	&\lesssim \frac{1}{t}|(x^{a}/r)H_{a}\left( [\varphi]_{-}\right)|+|\pt \left( [\varphi]_{-}\right)|\\
	&\lesssim  \frac{1}{t}\big|\big[(x^{a}/r)H_{a}\varphi\big]_{-}\big|+|\pt \left( [\varphi]_{-}\right)|\lesssim \frac{1}{t+r}\left|\Gamma \varphi\right|+|G|.
	\end{aligned}
	\end{equation}
	We see that~\eqref{est:ptvarphi-} follows from~\eqref{est:GammaGammahat} and~\eqref{est:pavarphi}.
	\end{proof}

\subsection{Estimates for the 3D linear wave equation}
In this subsection, we review several technical estimates for the 3D linear wave equation. First, we recall the following two $L^{\infty}$--$L^{\infty}$ estimates for the 3D linear wave equation from~\cite[Proposition 3.1]{AlinhacIndiana} and~\cite[Lemma 3.4]{Kata} respectively.
\begin{proposition}[\cite{AlinhacIndiana}]\label{pro:Linfinity}
	Let $u$ be the solution to the Cauchy problem
	\begin{equation*}
	-\Box u=G\quad \mbox{with}\quad (u,\pt u)_{|t=2}=(0,0),
	\end{equation*}
	and assume that $u$ is supported in $\mathcal{K}$. Assume the estimates 
	\begin{equation*}
	|G|\le C_{G} (t-r)^{-\mu}(t+r)^{-\nu},\quad \mu,\nu\ge 0.
	\end{equation*}
	Define $F_{\mu}(s)=1,\log s, s^{1-\mu}/(1-\mu)$ according to $\mu >1,=1,<1$, respectively, then the following pointwise estimates are true.
	\begin{enumerate}
		
		\item  For $\nu<2$, we have 
		\begin{equation*}
		|u(t,x)|\lesssim C_{G}(t+r)^{-\nu+1}F_{\mu}(t-r).
		\end{equation*}
		
		\item  For $\nu=2$, we have 
		\begin{equation*}
		|u(t,x)|\lesssim C_{G}(t+r)^{-1}\log (t+r)F_{\mu}(t-r).
		\end{equation*}
		
		\item  For $\nu >2$, we have 
		\begin{equation*}
		|u(t,x)|\lesssim C_{G}(t-r)^{-(\nu-2)}(t+r)^{-1}F_{\mu}(t-r).
		\end{equation*}
		
\end{enumerate}
\end{proposition}

\begin{proposition}[\cite{Kata,Kubota-Yokoyama}]\label{pro:Linwave}
	Let $u$ be the solution to the Cauchy problem
	\begin{equation*}
	-\Box u=G\quad \mbox{with}\quad (u,\pt u)_{|t=2}=(0,0).
	\end{equation*}
	Then we have 
	\begin{equation*}
	\langle t+r\rangle\langle t-r\rangle^{\frac{1}{2}}|u(t,x)|\lesssim \sup_{\tau \in [2,t]}\sup_{y\in \R^{3}}\langle \tau +|y|\rangle^{\frac{7}{2}}|G(\tau,y)|.
	\end{equation*}
\end{proposition}

\begin{remark}
	Since we lack the scaling vector field when treating massive Dirac or Klein-Gordon equations, the classical Klainerman-Sobolev inequality cannot be directly applied, and thus in many cases, we cannot get robust pointwise decay using weighted $L^2$ type bounds (see for instance \eqref{est:SobolevGlobal}). The important part of Proposition \ref{pro:Linwave} is that we do not rely on the scaling vector field, and it allows us to obtain sharp pointwise decay of $u$ in the sense that the difference of the homogeneity between $u$ decay and $G$ decay is $2$. We also note that Proposition \ref{pro:Linwave} is just a special version of \cite[Lemma 3.4]{Kata} with $\rho=0, \kappa={3/2}, \mu={1/2}$ therein.
\end{remark}

Second, we introduce the following standard energy $\mathcal{E}$ and conformal energy $\mathcal{F}$ for the 3D linear wave equation,
\begin{equation*}
\begin{aligned}
\mathcal{E}(t,u)&=\int_{\R^{3}}\left( (\pt u)^{2}+|\nabla_{x} u|^{2}\right)(t,x)\d x,\\
\mathcal{F}(t,u)&=\int_{\R^{3}}\left( u^{2}+(Su)^{2}+|\Omega u|^{2}+|Hu|^{2}\right)(t,x)\d x.
\end{aligned}
\end{equation*}

\begin{lemma}[\cite{AlinhacBook}]\label{le:waveenergy}
		Let $u$ be the solution to the Cauchy problem
	\begin{equation*}
	-\Box u=G\quad \mbox{with}\quad (u,\pt u)_{|t=2}=(u_{0},u_{1}).
	\end{equation*}
	Then the following estimates hold.
	\begin{enumerate}
		\item \emph{Standard energy estimate.} We have 
		\begin{equation*}
		\mathcal{E}(t,u)^{\frac{1}{2}}\lesssim \mathcal{E}(2,u)^{\frac{1}{2}}+\int_{2}^{t}\|G(\tau,x)\|_{L_{x}^{2}}\d \tau.
		\end{equation*}
		
		\item \emph{Conformal energy estimate.} We have 
		\begin{equation*}
	\mathcal{F}(t,u)^{\frac{1}{2}}\lesssim \mathcal{F}(2,u)^{\frac{1}{2}}+\int_{2}^{t}\left\|\langle \tau+r\rangle G(\tau,x)\right\|_{L_{x}^{2}}\d \tau.
		\end{equation*}
	\end{enumerate}
\end{lemma}

\begin{proof}
	See~\cite[Theorem 6.3]{AlinhacBook} and~\cite[Theorem 6.11]{AlinhacBook}.
	\end{proof}

Last, we introduce the following weighted conformal energy estimate for 3D linear wave equation. This result is motivated by a $(t-r)$-weighted basic energy estimate for the wave equation in \cite{Alinhac01b, Dong2006}. The proof is inspired by the previous works~\cite[Theorem 6.11]{AlinhacBook} and~\cite[Lemma 6.3.5]{HomanderBook}.
\begin{lemma}\label{le:weightcon}
	Let $u$ be the solution to the Cauchy problem
	\begin{equation*}
	-\Box u=G\quad \mbox{with}\quad (u,\pt u)_{|t=2}=(u_{0},u_{1}).
	\end{equation*}
	Assume that $u$ and $G$ are supported in $\mathcal{K}$, then we have 
	\begin{equation}\label{est:Con}
	\begin{aligned}
	\sum_{\Lambda\in \left\{{\rm{Id}},S,\Omega,H\right\}}\left\|\frac{\Lambda u}{(t-r)^{\frac{1}{2}}}\right\|_{L_{x}^{2}}\lesssim \left\| u_{0}\right\|_{H^{1}}+\|u_{1}\|_{L_{x}^{2}}+\int_{2}^{t}\left\|\frac{\tau G(\tau,x)}{(\tau-r)^{\frac{1}{2}}}\right\|_{L_{x}^{2}}\d \tau.
	\end{aligned}
	\end{equation}
\end{lemma}

\begin{remark}
We first comment that the weight $(t-r)^{-1/2}$ is not typical, as the conformal energy estimate \eqref{est:Con} is also true if we replace $(t-r)^{-1/2}$ by $(t-r)^{-\gamma}$ with $\gamma > 0$ on both sides of \eqref{est:Con}. Second, the advantage of \eqref{est:Con} is that it allows us to utilise the $(t-r)/t$ decay in the source term $G$, and in turn this $(t-r)/t$ decay can be inherited by $\Lambda u$ with $\Lambda\in \left\{{\rm{Id}},S,\Omega,H\right\}$.
\end{remark}

\begin{proof}[Proof of Lemma~\ref{le:weightcon}]
	Consider the following weighted conformal energy, 
	\begin{equation*}
	\begin{aligned}
	\mathcal{E}_{con}(t,u)
	&=\frac{1}{2}\int_{\R^{3}}\left[(r^{2}+t^{2})|\partial u|^{2}+4rt\pt u\partial_{r}u\right](t-r)^{-1}\d x\\
	&+\int_{\R^{3}}\left[2tu \pt u-u^{2}+{\rm{div}}(xu^{2})\right](t-r)^{-1}\d x.
	\end{aligned}
	\end{equation*}
	
	\textbf{Step 1.} Identity for $\mathcal{E}_{con}$. We claim that 
	\begin{equation}\label{equ:Con}
	\mathcal{E}_{con}(t,u)=\frac{1}{2}\int_{\R^{3}}\left[(Su+2u)^{2}+|Hu|^{2}+|\Omega u|^{2}\right](t-r)^{-1}\d x.
	\end{equation}
	Indeed, from the definition of $S$, $H$ and $\Omega$, we have 
	\begin{equation*}
	\begin{aligned}
	&\frac{1}{2}\left[(r^{2}+t^{2})|\partial u|^{2}+4rt(\pt u)(\partial_{r}u)\right]\\
	=&\frac{1}{2}(r^{2}+t^{2})\left[(\pt u)^{2}+(\partial_{r} u)^{2}+\left|\frac{\Omega u}{r}\right|^{2}\right]+2rt (\pt u)(\partial_{r} u)\\
	=&\frac{1}{2}\left[(Su)^{2}+(t\partial_{r}u+r\pt u)^{2}\right]+\frac{1}{2}(r^{2}+t^{2})\left|\frac{\Omega u}{r}\right|^{2}=\frac{1}{2}\left[(Su)^{2}+|Hu|^{2}+|\Omega u|^{2}\right].
	\end{aligned}
	\end{equation*}
	Then, by an elementary computation, we check that 
	\begin{equation*}
	2tu \pt u-u^{2}+{\rm{div}}(xu^{2})=2uSu+2u^{2}.
	\end{equation*}
	Combining the above two identities and the definition of $\mathcal{E}_{con}$, we obtain~\eqref{equ:Con}.
	
	\smallskip
	\textbf{Step 2.} First estimate for $\mathcal{E}_{con}$. We claim that 
	\begin{equation}\label{est:Con1}
	\mathcal{E}_{con}^{\frac{1}{2}}(t,u)-\mathcal{E}_{con}^{\frac{1}{2}}(2,u)\lesssim \int_{2}^{t}\left\|\frac{\tau G(\tau,x)}{(\tau-r)^{\frac{1}{2}}}\right\|_{L_{x}^{2}}\d \tau. 
	\end{equation}
	Indeed, by a tedious computation (see more details in~\cite[Page 99]{AlinhacBook}), we have 
	\begin{equation*}
	\begin{aligned}
	-\Box u\left(K_{0}u+2tu\right)
	&=\frac{1}{2}\pt \left[(r^{2}+t^{2})|\partial u|^{2}+4rt(\pt u)(\partial_{r} u)\right]\\
	&+\pt \left[2tu\pt u-u^{2}+{\rm{div}}(xu^{2})\right]-{\rm{div}}\left[2tu \nabla u+\pt (xu^{2})\right]\\
	&+{\rm{div}}\left[tx\left(-2(\pt u)^{2}+|\partial u|^{2}\right)-(r^{2}+t^{2})(\pt u)\nabla u-2rt(\partial_{r}u)\nabla u\right],
	\end{aligned}
	\end{equation*}
	where $K_{0}=(r^{2}+t^{2})\pt +2rt\partial_{r}$. Based on the above identity, we have 
	\begin{equation*}
	\begin{aligned}
	&-\Box u\left(K_{0}u+2tu\right)(t-r)^{-1}\\
	&=\mathcal{I}_{1}+\mathcal{I}_{2}+\mathcal{I}_{3}+\frac{1}{2}\pt \left\{\left[(r^{2}+t^{2})|\partial u|^{2}+4rt(\pt u)(\partial_{r} u)\right](t-r)^{-1}\right\}\\
		&+\pt \left\{\left[2tu\pt u-u^{2}+{\rm{div}}(xu^{2})\right](t-r)^{-1}\right\}-{\rm{div}}\left\{\left[2tu \nabla u+\pt (xu^{2})\right](t-r)^{-1}\right\}\\
	&+{\rm{div}}\left\{\left[tx\left(-2(\pt u)^{2}+|\partial u|^{2}\right)-(r^{2}+t^{2})(\pt u)\nabla u-2rt(\partial_{r}u)\nabla u\right](t-r)^{-1}\right\},
	\end{aligned}
	\end{equation*}
	where
	\begin{equation*}
	\begin{aligned}
	\mathcal{I}_{1}&=-\frac{1}{2}\left[(r^{2}+t^{2})|\partial u|^{2}+4rt(\pt u)(\partial_{r} u)\right] \pt \left((t-r)^{-1}\right),\\
	\mathcal{I}_{2}&=\left[2tu \nabla u+\pt (xu^{2})\right]\cdot \nabla\left((t-r)^{-1}\right)-\left[2tu\pt u-u^{2}+{\rm{div}}(xu^{2})\right] \pt \left((t-r)^{-1}\right),\\
	\mathcal{I}_{3}&=-\left[tx\left(-2(\pt u)^{2}+|\partial u|^{2}\right)-(r^{2}+t^{2})(\pt u)\nabla u-2rt(\partial_{r}u)\nabla u\right]\cdot \nabla\left((t-r)^{-1}\right).
	\end{aligned}
	\end{equation*}
	Integrating the above identity over $\R^{3}$, we have 
	\begin{equation}\label{est:dtEcon}
	\frac{\d}{\d t}{\mathcal{E}}_{con}+\int_{\R^{3}}(\mathcal{I}_{1}+\mathcal{I}_{2}+\mathcal{I}_{3})\d x 
	=\int_{\R^{3}}\left[G(t, x)\left(K_{0}u+2tu\right)(t-r)^{-1}\right]\d x.
	\end{equation}
	Note that, the following identities are true:
	\begin{equation}\label{equ:ptt-r}
	\pt ((t-r)^{-1})=- (t-r)^{-2}\quad \mbox{and}\quad 
	\nabla_{x}((t-r)^{-1})= (t-r)^{-2}\frac{x}{r}.
	\end{equation}
	Therefore, by an elementary computation, we obtain
	\begin{equation*}
	\begin{aligned}
	\mathcal{I}_{1}&=\frac{1}{2}\left[(r^{2}+t^{2})|\partial u|^{2}+4rt(\pt u)(\partial_{r} u)\right](t-r)^{-2},\\
	\mathcal{I}_{2}&=2\left[u^{2}+(t+r)u(\partial_{t}u+\partial_{r} u)\right](t-r)^{-2},\\
	\mathcal{I}_{3}&=\left[tr(2(\pt u)^{2}+2(\partial_{r}u)^{2}-|\partial u|^{2})+(r^{2}+t^{2})(\pt u)(\partial_{r}u)\right](t-r)^{-2}.
	\end{aligned}
	\end{equation*}
	It follows from $|\partial u|^{2}=(\pt u)^{2}+(\partial_{r} u)^{2}+r^{-2}|\Omega u|^{2}$ that
	\begin{equation*}
	\begin{aligned}
	&\int_{\R^{3}}\left(\mathcal{I}_{1}+\mathcal{I}_{2}+\mathcal{I}_{3}\right)\d x\\
	&=\frac{1}{2}\int_{\R^{3}}(t-r)^{-2}\left[\left((t+r)(\pt u+\partial_{r}u)+2u\right)^{2}+\big(1-\frac{t}{r}\big)^{2}\left|\Omega u\right|^{2}\right]\d x\ge 0.
	\end{aligned}
	\end{equation*}
	From~\eqref{est:dtEcon} and the above inequality, we know that 
	\begin{equation*}
	\frac{\d }{\d t}\mathcal{E}_{con}\le \int_{\R^{3}}\left|G\left(K_{0}u+2tu\right)(t-r)^{-1}\right|\d x.
	\end{equation*}
	On the other hand, from the definition of $K_{0}=(r^{2}+t^{2})\pt +2rt\partial_{r}$, we find
	\begin{equation*}
	\begin{aligned}
	\left(K_{0}u+2tu\right)^{2}
	&=\left[t\left(Su+u\right)+r\left(r\pt u+t\partial_{r}u\right)\right]^{2}\\
	&\le 2(t^{2}+r^{2}) \left[(Su+2u)^{2}+(t\partial_{r}u+r\pt u)^{2}\right]\\
	&\le 2(t^{2}+r^{2})\left[(Su+2u)^{2}+|Hu|^{2}\right]\le 4t^{2}\left[(Su+2u)^{2}+|Hu|^{2}\right],
	\end{aligned}
	\end{equation*}
	which implies
	\begin{equation*}
	\frac{\d }{\d t}\mathcal{E}_{con}(t,u)\le 2\left\|\frac{tG(t,x)}{(t-r)^{\frac{1}{2}}}\right\|_{L_{x}^{2}}\mathcal{E}_{con}^{\frac{1}{2}}(t,u).
	\end{equation*}
	Combining the above inequality with the Cauchy-Schwarz inequality, we obtain~\eqref{est:Con1}.
	
	\smallskip
	\textbf{Step 3.} Second estimate for $\mathcal{E}_{con}$. We claim that 
	\begin{equation}\label{est:Con2}
	\int_{\R^{3}}\frac{u^{2}}{( t-r)}\d x \le \mathcal{E}_{con}(t,u).
	\end{equation}
	Indeed, on one hand, by an elementary computation, we check that 
	\begin{equation*}
	\begin{aligned}
	&\frac{1}{2}(r^{2}+t^{2})|\partial u|^{2}+2rt(\pt u)(\partial_{r}u)+2tu\pt u-u^{2}\\
	&=\frac{1}{2}(r^{2}+t^{2})\left[\pt u+\frac{2t}{r^{2}+t^{2}}(u+r\partial_{r}u)\right]^{2}+\frac{1}{2r^{2}}\left(r^{2}+t^{2}\right)\left|\Omega u\right|^{2}\\
	&+\frac{1}{2}(r^{2}+t^{2})\left(\partial_{r} u\right)^{2}-\frac{2t^{2}}{r^{2}+t^{2}}(u+r\partial_{r}u)^{2}-u^{2}.
	\end{aligned}
	\end{equation*}
	Setting $v=ru$, we can rewrite the last three terms above as:
	\begin{equation*}
	\begin{aligned}
	&\frac{1}{2}(r^{2}+t^{2})\left(\partial_{r} u\right)^{2}-\frac{2t^{2}}{r^{2}+t^{2}}(u+r\partial_{r}u)^{2}-u^{2}\\
	=&\frac{(r^{2}-t^{2})^{2}}{2r^{2}(r^{2}+t^{2})}(\partial_{r}v)^{2}
	+\frac{1}{2r^{4}}\left(t^{2}-r^{2}\right)v^{2}-\frac{r^{2}+t^{2}}{r^{3}}v\partial_{r}v.
	\end{aligned}
	\end{equation*}
	Therefore, from the polar coordinate transformation $\d x=r^{2}\d r \d \omega$, we have 
	\begin{equation}\label{est:Con211}
	\begin{aligned}
	&\int_{\R^{3}}\left[\frac{1}{2}(r^{2}+t^{2})|\partial u|^{2}+2rt(\pt u)(\partial_{r}u)+2tu\pt u-u^{2}\right](t-r)^{-1}\d x\\
	&\ge \int_{\mathbb{S}^{2}}\int_{0}^{\infty}\left[\frac{1}{2r^{2}}(t^{2}-r^{2})v^{2}-\frac{r^{2}+t^{2}}{r}v\partial_{r}v\right](t-r)^{-1}\d r\d \omega.
	\end{aligned}
	\end{equation}
	Note that $v=O(r)$ near $r=0$, and so we can integrate by parts the last term with respect to $r$ between $r=0$ and $r=t-1$:
	\begin{equation}\label{equ:Con211}
	\begin{aligned}
	&\int_{{\mathbb{S}}^{2}}\int_{0}^{t-1}\left(-\frac{r^{2}+t^{2}}{r}v\partial_{r}v\right)(t-r)^{-1}\d r \d \omega\\
	&=\int_{{\mathbb{S}}^{2}}\int_{0}^{t-1}v^{2}\partial_{r}\left[\frac{r^{2}+t^{2}}{2r}(t-r)^{-1}\right]\d r \d \omega\\
	&=\int_{{\mathbb{S}}^{2}}\int_{0}^{t-1}v^{2}\left(\frac{1}{2r^{2}}(r^{2}-t^{2})+\frac{1}{2}\frac{r^{2}+t^{2}}{r(t-r)}\right)(t-r)^{-1}\d r \d \omega.
	\end{aligned}
	\end{equation}
	Using~\eqref{est:Con211} and~\eqref{equ:Con211}, we infer that 
	\begin{equation}\label{est:Con21}
	\begin{aligned}
	&\int_{\R^{3}}\left[\frac{1}{2}(r^{2}+t^{2})|\partial u|^{2}+2rt(\pt u)(\partial_{r}u)+2tu\pt u-u^{2}\right](t-r)^{-1}\d x\\
	&\ge \frac{1}{2}\int_{{\mathbb{S}}^{2}}\int_{0}^{t-1}\frac{r^{2}+t^{2}}{r(t-r)^{2}}v^{2}\d r \d \omega\ge \frac{1}{2}\int_{\R^{3}}\frac{r^{2}+t^{2}}{r(t-r)^{2}}u^{2}\d x.
	\end{aligned}
	\end{equation}
	On the other hand, using again the integration by parts and~\eqref{equ:ptt-r}, we have 
	\begin{equation}\label{equ:Condivxu2}
	\begin{aligned}
	&\int_{\R^{3}}\left[{\rm{div}}(xu^{2})\right](t-r)^{-1}\d x\\
	&=-\int_{\R^{3}}(xu^{2})\cdot \nabla (t-r)^{-1}\d x 
	=-\frac{1}{2}\int_{\R^{3}}\frac{2ru^{2}}{(t-r)^{2}}\d x.
	\end{aligned}
	\end{equation}
	Combining~\eqref{est:Con211} and~\eqref{equ:Condivxu2} with ${\rm{supp}}u\subseteq \mathcal{K}$, we obtain
	\begin{equation*}
	\mathcal{E}_{con}(t,u)\ge \frac{1}{2}\int_{\R^{3}}\frac{(t+r)u^{2}}{r(t-r)}\d x\ge \int_{\R^{3}}\frac{u^{2}}{(t-r)}\d x,
	\end{equation*}
	which means~\eqref{est:Con2}.
	
	\smallskip 
	\textbf{Step 4.} Conclusion. From~\eqref{equ:Con},~\eqref{est:Con2} and the triangle inequality, we obtain
	\begin{equation}\label{est:Con3}
	\begin{aligned}
	&\int_{\R^{3}}\frac{|Su|^{2}}{(t-r)}\d x+\int_{\R^{3}}\frac{u^{2}}{(t-r)}\d x+\int_{\R^{3}}\frac{|H u|^{2}}{(t-r)}\d x+\int_{\R^{3}}\frac{|\Omega u|^{2}}{(t-r)}\d x\\
	&\lesssim \int_{\R^{3}}\big[\left(Su+2u\right)^{2}+u^{2}+|Hu|^{2}+\left|\Omega u\right|^{2}\big](t-r)^{-1}\d x \lesssim \mathcal{E}_{con}(t,u).
	\end{aligned}
	\end{equation}
	We see that~\eqref{est:Con} follows from~\eqref{est:Con1} and~\eqref{est:Con3}.
	\end{proof}

\subsection{Estimates for the 3D linear Klein-Gordon equation}
In this subsection, we review several technical estimates of 3D linear Klein-Gordon equation for further reference. First, we recall the following pointwise estimate from~\cite[Theorem 1]{Geor}.

\begin{proposition}[\cite{Geor}]\label{prop:decayKG}
		Let $u$ be the solution to the Cauchy problem
	\begin{equation*}
	-\Box u+u=G\quad \mbox{with}\quad (u,\pt u)_{|t=2}=(u_{0},u_{1}).
	\end{equation*}
	Then we have 
	\begin{equation*}
	\begin{aligned}
	\langle t+r\rangle^{\frac{3}{2}}|u(t,x)|&\lesssim 	\sum_{k=0}^{\infty}\sum_{|I|\le 5}\left\|\langle |x|\rangle^{\frac{3}{2}} \zeta_{k}(|x|)Z^{I}u(0,x)\right\|_{L_{x}^{2}}\\
	&+\sum_{k=0}^{\infty}\sum_{|I|\le 4}\max_{2\le \tau\le t}\zeta_{k}(\tau)\left\|( \tau+|x|)Z^{I}G(\tau,x)\right\|_{L_{x}^{2}}.
	\end{aligned}
	\end{equation*}
\end{proposition}

Second, from the above Proposition and the Cauchy-Schwarz inequality, we have the following simplified version of Proposition~\ref{prop:decayKG}.

\begin{corollary}\label{co:Klein}
	In the context of Proposition~\ref{prop:decayKG}, let $\delta_{0}>0$ and assume
	\begin{equation}\label{est:assumg}
	\sum_{|I|\le 4}\max_{2\le \tau\le t}\tau^{\delta_{0}}\left\|( \tau+|x|)Z^{I}G(\tau,x)\right\|_{L_{x}^{2}}\le C_{G}.
	\end{equation}
	Then we have 
		\begin{equation*}
	\langle t+r\rangle^{\frac{3}{2}} |u(t,x)|\lesssim \frac{C_{G}}{1-2^{-\delta_{0}}}+\sum_{|I|\le 5}\left\|\langle |x|\rangle^{\frac{3}{2}} \log (2+|x|)Z^{I}u(0,x)\right\|_{L_{x}^{2}}.
	\end{equation*}
\end{corollary}

\begin{proof}
	We refer to~\cite[Corollary 2.12]{DLMY} for the details of the proof.
	\end{proof}

Next, we introduce the extra decay for the 3D linear Klein-Gordon equation, which was first observed by Klainerman \cite{Klainerman93}.

\begin{lemma}[\cite{Klainerman93}]\label{le:decayKG}
		Let $u$ be the solution to the Cauchy problem
	\begin{equation*}
	-\Box u+u=G\quad \mbox{with}\quad (u,\pt u)_{|t=2}=(u_{0},u_{1}).
	\end{equation*}
	Then we have 
	\begin{equation*}
	|u|\lesssim \frac{\langle t-r\rangle}{\langle t+r\rangle}\left(|\partial Zu|+|\partial u|\right)+|G|.
	\end{equation*}
\end{lemma}

\begin{proof}
	By an elementary computation, we rewrite the d'Alembert operator $-\Box$ as 
	\begin{equation*}
	-\Box=\frac{(t+r)(t-r)}{t^{2}}\partial_{0} \partial_{0}+\frac{3}{t}\partial_{0} +\frac{x^{a}}{t^{2}}\partial_{0} H_{a} -\frac{1}{t}\partial^{a}H_{a}-\frac{x^{a}}{t^{2}}\partial_{a}.
	\end{equation*}
	Therefore, from the equation of $u$ and the triangle inequality, we obtain
	\begin{align*}
	|u|\lesssim |\Box u|+|G|\lesssim \frac{\langle t-r\rangle}{\langle t+r\rangle}\left(\left|\partial Zu\right|+\left|\partial u\right|\right)+|G|,\quad \mbox{for}\ r\le 3t.
	\end{align*}
	On the other hand, using again the equation of $u$, we obtain
	\begin{equation*}
		|u|\lesssim |\Box u|+|G|\lesssim \frac{\langle t-r\rangle}{\langle t+r\rangle}\left|\partial \partial u\right|+|G|,\quad \mbox{for}\ r\ge 3t.
	\end{equation*}
	Combining the above two inequalities, we complete the proof of Lemma~\ref{le:decayKG}.
	\end{proof}

\section{Proof of Theorem~\ref{thm:massless}}\label{S:thmmassless}

In this section, we prove the existence of global-in-time solution $(\psi,\phi)$ of~\eqref{equ:3DDKG} satisfying~\eqref{est:theorem1point}-\eqref{est:theorem1sca} in Theorem~\ref{thm:massless} via a bootstrap argument. 

\subsection{Bootstrap assumption}
 Let $N\in \mathbb{N}^{+}$ with $N\ge 7$ and $0<\delta\ll 1$. The proof of Theorem~\ref{thm:massless} relies on a bootstrap argument of high-order energy and pointwise decay of solution $(\psi,\phi)$. More precisely, we consider the following bootstrap assumption: for $C\gg 1$ and $0<\epsilon\ll C^{-1}$ to be chosen later, 
\begin{equation}\label{est:Boot1}
\left\{\begin{aligned}
\sum_{|I|\le N}\left(\mathcal{E}_{D}^{\delta}\big(t,\widehat{\Gamma}^{I}\psi\big)^{\frac{1}{2}}+\left\|\partial\Gamma^{I}\phi\right\|_{L_{x}^{2}}\right)&\le C\epsilon,\\
\sum_{|I|\le N-4}(t+r)^{1-\delta}(t-r)^{\frac{1}{2}+4\delta}\big|\widehat{\Gamma}^{I}\psi\big|&\le C\epsilon,\\
\sum_{|I|\le N-3}(t+r)(t-r)^{2\delta}\left|\Gamma^{I}\phi\right|&\le C\epsilon,\\
\sum_{|I|\le N-4}(t+r)(t-r)^{6\delta}\left|\Gamma^{I}\phi\right|&\le C\epsilon.
\end{aligned}\right.
\end{equation}

For all initial data $(\psi_{0},\vec{\phi}_{0})$ satisfying~\eqref{est:small}, we set 
\begin{equation}\label{def:T*}
T_{*}:=T_{*}(\psi_{0},\vec{\phi}_{0})=\sup \left\{t\in [2,\infty):(\psi,\phi) \ \mbox{satisfy}~\eqref{est:Boot1}\ \mbox{on}\ [2,t)\right\}>2.
\end{equation}

The following proposition is the main part of the proof of Theorem~\ref{thm:massless}.
\begin{proposition}\label{pro:main1}
	For any compactly supported initial data $(\psi_{0},\vec{\phi}_{0})$ satisfying the smallness condition~\eqref{est:small} in Theorem~\ref{thm:massless}, we have $T_{*}(\psi_{0},\vec{\phi}_{0})=\infty$.
\end{proposition}

The rest of the section is organized as follows. First, in \S\ref{SS:Keypoint}, we introduce several key estimates for $(\psi,\phi)$. Later, in \S\ref{SS:EndThm1}, we prove Proposition~\ref{pro:main1} and then finish the proof of Theorem~\ref{thm:massless} from Proposition~\ref{pro:main1}. In what follows, the implied constants in $\lesssim$ do not depend on the
constants $C$ and $\epsilon$ appearing in the bootstrap assumption~\eqref{est:Boot1}.

\subsection{Key estimates for $(\psi,\phi)$}\label{SS:Keypoint}

First of all, from~\eqref{equ:3DDKG} and~\eqref{equ:comm}, for any $I\in \mathbb{N}^{11}$ with $|I|\le N$, there exist constants $\left\{a_{J}^{I}\right\}_{|J|\le |I|}$ and $\left\{b_{J}^{I}\right\}_{|J|\le |I|}$ (independent of $(\psi,\phi)$) such that
\begin{equation}\label{equ:3DDKGuse}
\begin{aligned}
&-i\gamma^{\mu}\partial_{\mu}\widehat{\Gamma}^{I}\psi=a_{J}^{I}\widehat{\Gamma}^{J}\left(\phi\psi\right)\quad \mbox{and}\quad -\Box \Gamma^{I}\phi=b_{J}^{I}\widehat{\Gamma}^{J}\left(\psi^{*}\gamma^{0}\psi\right).
\end{aligned}
\end{equation}
The above identities will be frequently used in this section.

Note that, from~\eqref{est:GammaGammahat},~\eqref{est:Klainer} and~\eqref{est:Boot1}, we directly have 
\begin{equation}\label{est:Firstpsi}
\sum_{|I|\le N-2}\big|\widehat{\Gamma}^{I}\psi\big|\lesssim C\epsilon(t+r)^{-1}(t-r)^{-\frac{1}{2}}.
\end{equation}

Note also that, from~\eqref{est:Boot1} and the Cauchy-Schwarz inequality, we have 
\begin{equation}\label{est:Ghostpsi}
\sum_{|I|\le N}\int_{2}^{t}\tau^{-\frac{1}{2}-\delta}\left\|\frac{\big[\widehat{\Gamma}^{I}\psi\big]_{-}}{(\tau-r)^{\frac{1}{2}+\delta}}\right\|_{L_{x}^{2}}\d \tau\lesssim C\epsilon.
\end{equation}

Before proceeding further, let us introduce an auxiliary variable $\Psi$, which is the solution to a 3D linear wave equation with a specific source term. More precisely, following~\cite[Section 2]{Boura00}, for any solution $(\psi,\phi)$ of the 3D massless Dirac--Klein-Gordon system~\eqref{equ:3DDKG}, we denote by $\Psi=\Psi(t,x)$ the solution to the following 3D wave equation, 
\begin{equation}\label{equ:defPsi}
\left\{\begin{aligned}
-\Box\Psi(t,x)&=i\gamma^{\mu}\partial_{\mu}\psi(t,x),\quad \mbox{for}\ (t,x)\in [2,\infty)\times \R^{3},\\
(\Psi,\pt \Psi)_{|t=2}&=(0,i\gamma^{0}\psi),\quad \quad \ \mbox{for}\ x\in \R^{3}.
\end{aligned}\right.
\end{equation}
Then $-i\gamma^{\mu}\partial_{\mu}\Psi=\psi$, since $-i\gamma^{\nu}\partial_{\nu}\left(i\gamma^{\mu}\partial_{\mu}\Psi+\psi\right)=0$ and $\left(i\gamma^{\mu}\partial_{\mu}\Psi+\psi\right)_{|t=2}=0$. Recall that, the following lemma is true for 3D Dirac--Klein-Gordon system~\eqref{equ:3DDKG}.
\begin{lemma}
	For any $I\in \mathbb{N}^{11}$ with $|I|\le N-2$, we have 
	\begin{equation}\label{equ:Hidden}
	\big|\big[\widehat{\Gamma}^{I}\psi\big]_{-}\big|\lesssim
	\sum_{|J|\le |I|}\sum_{a=1}^{3}\big|G_{a}{\widehat{\Gamma}}^{I}\Psi\big|.
	\end{equation}
\end{lemma}

\begin{proof}
	First, by an elementary computation, we have 
	\begin{equation*}
	-i\gamma^{\mu}\partial_{\mu}=-i\gamma^{b}G_{b}-i(\gamma^{0}-\gamma^{b}\omega_{b})\pt.
	\end{equation*}
Therefore, from~\eqref{equ:def+-},~\eqref{equ:comm} and $-i\gamma^{\mu}\partial_{\mu}\Psi=\psi$, there exist $\left\{c_{J}^{I}\right\}_{|J|\le |I|}$ and $\left\{d_{J}^{I}\right\}_{|J|\le |I|}$ (independent of $\psi$) such that
	\begin{equation*}
	\begin{aligned}
	\big[\widehat{\Gamma}^{I}\psi\big]_{-}=
	c_{J}^{I}\left(I_{4}-w_{a}\gamma^{0}\gamma^{a}\right)\gamma^{b}G_{b}\widehat{\Gamma}^{J}\Psi+d_{J}^{I}\left(I_{4}-w_{a}\gamma^{0}\gamma^{a}\right)\left(\gamma^{0}-\gamma^{b}\omega_{b}\right)\pt \widehat{\Gamma}^{J}\varphi.
	\end{aligned}
	\end{equation*}
	From~\eqref{equ:gamma}, we compute
	\begin{equation*}
	\begin{aligned}
	&\left(I_{4}-w_{a}\gamma^{0}\gamma^{a}\right)\left(\gamma^{0}-\gamma^{b}\omega_{b}\right)\\
	&=\gamma^{0}+\omega_{a}\omega_{b}\gamma^{0}\gamma^{a}\gamma^{b}=\bigg(1-\sum_{a=1}^{3}\omega_{a}^{2}\bigg)\gamma^{0}+\sum_{a<b}\omega_{a}\omega_{b}\gamma^{0}\left(\gamma^{a}\gamma^{b}+\gamma^{b}\gamma^{a}\right)=0.
	\end{aligned}
	\end{equation*}
	Combining the above two identities, we obtain
	\begin{equation}\label{equ:Hiddenclam}
	\big[\widehat{\Gamma}^{I}\psi\big]_{-}= c_{J}^{I}\left(I_{4}-w_{a}\gamma^{0}\gamma^{a}\right)\gamma^{b}G_{b}\widehat{\Gamma}^{J}\Psi.
	\end{equation}
	 We see that~\eqref{equ:Hidden} follows from~\eqref{equ:Hiddenclam} and the triangle inequality.
\end{proof}

Second, we deduce the following pointwise estimate for $\Gamma^{I}\Psi$.
\begin{lemma}
		For all $t\in [2,T_{*}(\psi_{0},\vec{\phi}_{0}))$, we have 
		\begin{equation}\label{est:Psi}
		\sum_{|I|\le N-3}\left|\Gamma^{I}\Psi\right|\lesssim  \left(\epsilon+C^{2}\epsilon^{2}\right) (t+r)^{-1+\delta}(t-r)^{\frac{1}{2}-6\delta}.
		\end{equation}
\end{lemma}

\begin{proof}
	We decompose the auxiliary variable $\Psi$ as,
	\begin{equation*}
	\Psi(t,x)=\Psi_{free}(t,x)+\Psi_{inho}(t,x)
	\end{equation*}
	where $\Psi_{free}$ and $\Psi_{inho}$ are the solutions for the following 3D linear homogeneous
	and inhomogeneous wave equations, respectively,
	\begin{equation}\label{equ:Psi1Psi2}
		\left\{\begin{aligned}
	-\Box \Psi_{free}&=0, \quad \quad \ \ \mbox{with}\ \ (\Psi_{free},\partial_{t}\Psi_{free})_{|t=2}=(0,i\gamma^{0}\psi),\\
-\Box \Psi_{inho}&=-\phi\psi,\quad  \mbox{with}\ \ (\Psi_{inho},\partial_{t}\Psi_{inho})_{|t=2}=(0,0).
	\end{aligned}\right.
	\end{equation}
	First, from~\eqref{est:small},~\eqref{equ:Psi1Psi2} and Lemma~\ref{le:waveenergy}, we have
	\begin{equation*}
	\sum_{|I|\le N-2}\left(\mathcal{E}(t,\Gamma^{I}\Psi_{free})^{\frac{1}{2}}+\mathcal{F}(t,\Gamma^{I}\Psi_{free})^{\frac{1}{2}}\right)\lesssim \|\psi_{0}\|_{H^{N}}\lesssim \epsilon,
	\end{equation*}
	which implies
	\begin{equation*}
	\sum_{|I|\le N-1}\left\|\Gamma^{I}\Psi_{free}\right\|_{L_{x}^{2}}\lesssim \sum_{|I|\le N-2}\left(\mathcal{E}(t,\Gamma^{I}\Psi_{free})^{\frac{1}{2}}+\mathcal{F}(t,\Gamma^{I}\Psi_{free})^{\frac{1}{2}}\right)\lesssim \epsilon.
	\end{equation*}
	Based on the above inequality and~\eqref{est:Klainer}, we obtain
	\begin{equation}\label{est:Psi1}
	\sum_{|I|\le N-3}\left|\Gamma^{I}\Psi_{free}\right|\lesssim \epsilon (t+r)^{-1}(t-r)^{-\frac{1}{2}}.
	\end{equation}
	Second, from the Leibniz rule,~\eqref{est:GammaGammahat},~\eqref{est:Boot1},~\eqref{est:Firstpsi} and $N\ge 7$, we have 
	\begin{equation*}
	\begin{aligned}
	\sum_{|I|\le N-3}\left|\Gamma^{I}(\phi\psi)\right|
	&\lesssim \sum_{\substack{|I_{1}|\le N-3\\ |I_{2}|\le N-4}}\left|\Gamma^{I_{1}}\phi\right|\big|\widehat{\Gamma}^{I_{2}}\psi\big|+ \sum_{\substack{|I_{1}|\le N-4\\ |I_{2}|\le N-3}}\left|\Gamma^{I_{1}}\phi\right|\big|\widehat{\Gamma}^{I_{2}}\psi\big|\\
	&\lesssim C^{2}\epsilon^{2}(t+r)^{-2+\delta}(t-r)^{-\frac{1}{2}-6\delta}.
	\end{aligned}
	\end{equation*}
	Therefore, from~\eqref{est:small},~\eqref{est:Klainer},~\eqref{equ:Psi1Psi2} and (i) of Proposition~\ref{pro:Linfinity}, one can obtain
	\begin{equation}\label{est:Psi2}
	\sum_{|I|\le N-3}\big|\Gamma^{I}\Psi_{inho}\big|
	\lesssim \epsilon (t+r)^{-1}(t-r)^{-\frac{1}{2}}+C^{2}\epsilon^{2}(t+r)^{-1+\delta}(t-r)^{\frac{1}{2}-6\delta}.
	\end{equation}
	We see that~\eqref{est:Psi} follows from~\eqref{est:Psi1} and~\eqref{est:Psi2}.
	\end{proof}

Third, we deduce the pointwise estimate for $\big[\widehat{\Gamma}^{I}\psi\big]_{-}$.

\begin{lemma}\label{le:Firstpsi}
	For all $t\in [2,T_{*}(\psi_{0},\vec{\phi}_{0}))$, the following estimates are true.
	\begin{enumerate}
		
		\item \emph{First estimate for $\big[\widehat{\Gamma}^{I}\psi\big]_{-}$.} We have 
		\begin{equation}\label{est:Firstpsi-}
		\sum_{|I|\le N-4}\big|\big[\widehat{\Gamma}^{I}\psi(t,x)\big]_{-}\big|\lesssim C\epsilon (t+r)^{-2+\delta}(t-r)^{\frac{1}{2}-6\delta}.
		\end{equation}
		
		\item \emph{Second estimate for $\big[\widehat{\Gamma}^{I}\psi\big]_{-}$.} We have 
		\begin{equation}\label{est:Secondpsi-}
		\sum_{|I|\le N-3}\big|\big[\widehat{\Gamma}^{I}\psi(t,x)\big]_{-}\big|\lesssim C\epsilon (t+r)^{-2+\delta}(t-r)^{\frac{1}{2}}.
		\end{equation}
		\end{enumerate}
\end{lemma}
\begin{proof}
	
	Proof of (i). Based on~\eqref{est:GammaGammahat},~\eqref{est:Ga},~\eqref{equ:Hidden} and~\eqref{est:Psi}, we obtain
	\begin{equation*}
	\begin{aligned}
	\sum_{|I|\le N-4}\big|\big[\widehat{\Gamma}^{I}\psi\big]_{-}\big|
	&\lesssim \sum_{|I|\le N-4}\sum_{a=1}^{3}\big|G_{a}\widehat{\Gamma}^{I}\Psi\big|\\
	&\lesssim (t+r)^{-1}\sum_{|I|\le N-4}\big|\Gamma\widehat{\Gamma}^{I}\Psi\big|\lesssim C\epsilon (t+r)^{-2+\delta}(t-r)^{\frac{1}{2}-6\delta},
	\end{aligned}
	\end{equation*}
	which means~\eqref{est:Firstpsi-}.
	
	Proof of (ii). Recall that, from~\eqref{est:GammaGammahat},~\eqref{est:Boot1},~\eqref{est:Firstpsi} and $N\ge 7$, we have 
	\begin{equation*}
	\begin{aligned}
	\sum_{|I|\le N-3}\big|\widehat{\Gamma}^{I}(\phi\psi)\big|
	&\lesssim C^{2}\epsilon^{2}(t+r)^{-2+\delta}(t-r)^{-\frac{1}{2}-6\delta}.
	\end{aligned}
	\end{equation*}
	
	Therefore, from~\eqref{est:ptvarphi-},~\eqref{est:Boot1},~\eqref{equ:3DDKGuse} and~\eqref{est:Firstpsi}, we obtain
	\begin{equation*}
	\begin{aligned}
	\sum_{|I|\le N-3}\big|\partial_{r} \big[\widehat{\Gamma}^{I}\psi\big]_{-}\big|
	&\lesssim (t+r)^{-1}\sum_{|I|\le N-2}\big|\widehat{\Gamma}^{I}\psi\big|+\sum_{|I|\le N-3} \big|\widehat{\Gamma}^{I}(\phi\psi)\big|\\
	&\lesssim C\epsilon (t+r)^{-2+\delta}(t-r)^{-\frac{1}{2}}.
	\end{aligned}
	\end{equation*}
	Denote $x=r\omega$ for $r=|x|$. Integrating the above inequality on $[r,t]$ and then using the fact that  ${\rm{supp}}\psi\subseteq \mathcal{K}$, we obtain
	\begin{equation*}
	\begin{aligned}
		\sum_{|I|\le N-3}\big|\big[\widehat{\Gamma}^{I}\psi\big]_{-}(t,x)\big|
		&\lesssim \sum_{|I|\le N-3}\int_{r}^{t}	\big|\partial_{\rho} \big[\widehat{\Gamma}^{I}\psi\big]_{-}(t,\rho\omega)\big|\d \rho\\
		&\lesssim C\epsilon\int_{r}^{t}(t+\rho)^{-2+\delta}(t-\rho)^{-\frac{1}{2}}\d \rho
		\lesssim C\epsilon (t+r)^{-2+\delta}(t-r)^{\frac{1}{2}},
	\end{aligned}
	\end{equation*}
	which means~\eqref{est:Secondpsi-}.
	\end{proof}

Last, we deduce the following weighted $L_{x}^{2}$ estimate for $\Gamma^{I}\phi$ for future reference.

\begin{lemma}
		For all $t\in [2,T_{*}(\psi_{0},\vec{\phi}_{0}))$, we have
		\begin{equation}\label{est:L2phi}
		\sum_{|I|\le N}\left\|\frac{\Gamma^{I}\phi}{(t-r)^{\frac{1}{2}}}\right\|_{L_{x}^{2}}\lesssim C\epsilon+C\epsilon t^{\delta}.
		\end{equation}
\end{lemma}

\begin{proof} The proof is based on an induction argument and we split it into four steps.
	
	\smallskip
	\textbf{Step 1.} Estimate for $\phi$. We claim that 
	\begin{equation}\label{est:L2phi0}
\left\|\frac{\phi}{(t-r)^{\frac{1}{2}}}\right\|_{L_{x}^{2}}+	\left\|\frac{\Gamma\phi}{(t-r)^{\frac{1}{2}}}\right\|_{L_{x}^{2}}\lesssim C\epsilon+C^{2}\epsilon^{2}t^{\delta}.
	\end{equation}
	Indeed, from~\eqref{est:Hidden},~\eqref{est:Boot1} and~\eqref{est:Firstpsi-}, for all $t\in [2,T_{*}(\psi_{0},\vec{\phi}_{0}))$, we have 
	\begin{equation*}
	\begin{aligned}
	\int_{2}^{t}\left\|\frac{\tau(\psi^{*}\gamma^{0}\psi)}{(\tau-r)^{\frac{1}{2}}}\right\|_{L_{x}^{2}}\d \tau
	&\lesssim \int_{2}^{t}\left\|\frac{\tau [\psi]_{-}}{(\tau-r)^{\frac{1}{2}}}\right\|_{L_{x}^{\infty}}\left\|\psi\right\|_{L_{x}^{2}}\d \tau\\
	&\lesssim C^{2}\epsilon^{2}\int_{2}^{t}\tau^{-1+\delta}\d \tau \lesssim C^{2}\epsilon^{2}t^{\delta}.
	\end{aligned}
	\end{equation*}
	It follows from~\eqref{est:small} and Lemma~\ref{le:weightcon} that
	\begin{equation}\label{est:L2phi01}
	\sum_{\Lambda\in\left\{{\rm{Id}},S,H,\Omega\right\}}\left\|\frac{\Lambda\phi}{(t-r)^{\frac{1}{2}}}\right\|_{L_{x}^{2}}\lesssim \epsilon+C^{2}\epsilon^{2}t^{\delta}.
	\end{equation}
	On the other hand, from~\eqref{est:Boot1} and ${\rm{supp}}\phi \subseteq \mathcal{K}$, we have
	\begin{equation}\label{est:L2phi02}
	\left\|\frac{\partial \phi}{(t-r)^{\frac{1}{2}}}\right\|_{L_{x}^{2}}\lesssim
	\left\|\partial \phi\right\|_{L_{x}^{2}}\lesssim C\epsilon.
	\end{equation}
	
	We see that~\eqref{est:L2phi0} follows from~\eqref{est:L2phi01} and~\eqref{est:L2phi02}.
	
	\smallskip
	\textbf{Step 2.} Key estimate for $\partial_{r} \big[\widehat{\Gamma}^{I}\psi\big]_{-}$. For $K\in \left\{1,\dots,N\right\}$, we claim that
	\begin{equation}\label{est:L2ppsi}
	\sum_{|I|\le K-1}\big\|\partial_{r} \big[\widehat{\Gamma}^{I}\psi\big]_{-}\big\|_{L_{x}^{2}}\lesssim C\epsilon t^{-1}+C\epsilon t^{-1}\sum_{|I|\le K-1}\left\|\frac{\Gamma^{I}\phi}{(t-r)^{\frac{1}{2}}}\right\|_{L_{x}^{2}}.
	\end{equation}
	Indeed, from $N\ge 7$,~\eqref{est:GammaGammahat},~\eqref{est:Boot1} and~\eqref{est:Firstpsi}, we have 
	\begin{equation*}
	\begin{aligned}
	\sum_{|I|\le K-1}\big\|\widehat{\Gamma}^{I}(\phi\psi)\big\|_{L_{x}^{2}}
		&\lesssim\sum_{\substack{|I_{1}|\le N-3\\ |I_{2}|\le K-1}}\left\|\Gamma^{I_{1}}\phi\right\|_{L_{x}^{\infty}}\big\|\widehat{\Gamma}^{I_{2}}\psi\big\|_{L_{x}^{2}}\\
	&+\sum_{\substack{|I_{1}|\le K-1\\|I_{2}|\le N-2}}\big\|(t-r)^{\frac{1}{2}}\widehat{\Gamma}^{I_{2}}\psi\big\|_{L_{x}^{\infty}} \left\|\frac{\Gamma^{I_{1}}\phi}{(t-r)^{\frac{1}{2}}}\right\|_{L_{x}^{2}}\\
	&\lesssim C^{2}\epsilon^{2} t^{-1}+C\epsilon t^{-1}\sum_{|I|\le K-1}\left\|\frac{\Gamma^{I}\phi}{(t-r)^{\frac{1}{2}}}\right\|_{L_{x}^{2}}.
	\end{aligned}
	\end{equation*}
	Therefore, from Lemma~\ref{le:ptvarphi},~\eqref{est:Boot1} and~\eqref{equ:3DDKGuse}, we obtain
	\begin{equation*}
	\begin{aligned}
		\sum_{|I|\le K-1}\big\|\partial_{r} \big[\widehat{\Gamma}^{I}\psi\big]_{-}\big\|_{L_{x}^{2}}
		&\lesssim t^{-1}\sum_{|I|\le K}\big\|\widehat{\Gamma}^{I}\psi\big\|_{L_{x}^{2}}+\sum_{|I|\le K-1}\big\|\widehat{\Gamma}^{I}(\phi\psi)\big\|_{L_{x}^{2}}\\
		&\lesssim C\epsilon t^{-1}+C\epsilon t^{-1}\sum_{|I|\le K-1}\left\|\frac{\Gamma^{I}\phi}{(t-r)^{\frac{1}{2}}}\right\|_{L_{x}^{2}},
		\end{aligned}
	\end{equation*}
	which means~\eqref{est:L2ppsi}.
	
	\smallskip
	\textbf{Step 3.} Key estimate for $\Gamma^{I}\phi$. For $K\in \left\{1,\dots,N\right\}$, we claim that 
	\begin{equation}\label{est:Keyphi}
	\sum_{|I|\le K}\left\|\frac{\Gamma^{I} \phi}{(t-r)^{\frac{1}{2}}}\right\|_{L_{x}^{2}}\lesssim C\epsilon t^{\delta}+C^{2}\epsilon^{2}\sum_{|I|\le K-1} \int_{2}^{t}\tau^{-1}\left\|\frac{\Gamma^{I}\phi}{(\tau-r)^{\frac{1}{2}}}\right\|_{L_{x}^{2}}\d \tau.
	\end{equation}
	Indeed, from $N\ge 7$,~\eqref{est:Hidden} and~\eqref{est:Gaga0}, we find
	\begin{equation*}
	\sum_{|I|\le K-1}\int_{2}^{t}\left\|\frac{\tau\Gamma^{I}(\psi^{*}\gamma^{0}\psi)}{(\tau-r)^{\frac{1}{2}}}\right\|_{L_{x}^{2}}\d \tau\lesssim \mathcal{G}_{1}+\mathcal{G}_{2},
	\end{equation*}
	where
	\begin{equation*}
	\begin{aligned}
		\mathcal{G}_{1}&=\sum_{|I_{1}|\le N-3}\sum_{|I_{2}|\le K-1}\int_{2}^{t}\left\|\frac{\tau \big[\widehat{\Gamma}^{I_{1}}\psi\big]_{-}\big|\widehat{\Gamma}^{I_{2}}\psi\big|}{(\tau-r)^{\frac{1}{2}}}\right\|_{L_{x}^{2}}\d \tau,\\
	\mathcal{G}_{2}&=\sum_{|I_{1}|\le K-1}\sum_{|I_{2}|\le N-2}\int_{2}^{t}\left\|\frac{\tau \big[\widehat{\Gamma}^{I_{1}}\psi\big]_{-}\big|\widehat{\Gamma}^{I_{2}}\psi\big|}{(\tau-r)^{\frac{1}{2}}}\right\|_{L_{x}^{2}}\d \tau.
	\end{aligned}
	\end{equation*}
	First, from~\eqref{est:Boot1} and~\eqref{est:Secondpsi-}, we have 
	\begin{equation*}
	\begin{aligned}
	\mathcal{G}_{1}
	&\lesssim \sum_{\substack{|I_{1}|\le N-3\\ |I_{2}|\le K-1}}\int_{2}^{t}\left\|\frac{\tau \big[\widehat{\Gamma}^{I_{1}}\psi\big]_{-}}{(\tau-r)^{\frac{1}{2}}}\right\|_{L_{x}^{\infty}}\big\|\widehat{\Gamma}^{I_{2}}\psi \big\|_{L_{x}^{2}} \d \tau\lesssim C^{2}\epsilon^{2}\int_{2}^{t}\tau^{-1+\delta}\d \tau \lesssim C^{2}\epsilon^{2} t^{\delta}.
	\end{aligned}
	\end{equation*}
	Second, from~\eqref{est:Hardy},~\eqref{est:Firstpsi} and~\eqref{est:L2ppsi}, we have
	\begin{equation*}
	\begin{aligned}
	\mathcal{G}_{2}
	&\lesssim \sum_{|I_{1}|\le K-1} \sum_{|I_{2}|\le N-2}
	\int_{2}^{t}\left\|\frac{ \big[\widehat{\Gamma}^{I_{1}}\psi\big]_{-}}{(\tau-r)}\right\|_{L_{x}^{2}}\left\|\tau (\tau-r)^{\frac{1}{2}} \big|\widehat{\Gamma}^{I_{2}}\psi\big|\right\|_{L_{x}^{\infty}}\d \tau\\
	&\lesssim C\epsilon \sum_{|I|\le K-1}\int_{2}^{t}\left\|\partial_{r} \big[\widehat{\Gamma}^{I}\psi\big]_{-}\right\|_{L_{x}^{2}}\d \tau\\
	&\lesssim C^{2}\epsilon^{2} \log t+C^{2}\epsilon^{2}\sum_{|I|\le K-1} \int_{2}^{t}\tau^{-1}\left\|\frac{\Gamma^{I}\phi}{(\tau-r)^{\frac{1}{2}}}\right\|_{L_{x}^{2}}\d \tau.
	\end{aligned}
	\end{equation*}
	Combining the above two inequalities, we have
	\begin{equation*}
	\sum_{|I|\le K-1}\int_{2}^{t}\left\|\frac{\tau\Gamma^{I}(\psi^{*}\gamma^{0}\psi)}{(\tau-r)^{\frac{1}{2}}}\right\|_{L_{x}^{2}}\d \tau\lesssim C^{2}\epsilon^{2} t^{\delta}+C^{2}\epsilon^{2}\sum_{|I|\le K-1} \int_{2}^{t}\tau^{-1}\left\|\frac{\Gamma^{I}\phi}{(\tau-r)^{\frac{1}{2}}}\right\|_{L_{x}^{2}}\d \tau.
	\end{equation*}
	Based on the above inequality,~\eqref{est:small},~\eqref{equ:3DDKGuse} and Lemma~\ref{le:weightcon}, one can obtain
	\begin{equation}\label{est:L21}
	\begin{aligned}
	&\sum_{\Lambda\in \left\{{\rm{Id}},S,H,\Omega\right\}}\sum_{|I|\le K-1}\left\|\frac{\Lambda \Gamma^{I}\phi}{(t-r)^{\frac{1}{2}}}\right\|_{L_{x}^{2}}\\
	&\lesssim \epsilon+C^{2}\epsilon^{2} t^{\delta}+C^{2}\epsilon^{2}\sum_{|I|\le K-1} \int_{2}^{t}\tau^{-1}\left\|\frac{\Gamma^{I}\phi}{(\tau-r)^{\frac{1}{2}}}\right\|_{L_{x}^{2}}\d \tau.
	\end{aligned}
	\end{equation}
	On the other hand, from~\eqref{est:Boot1} and ${\rm{supp}}\phi \subseteq \mathcal{K}$, one can obtain
	\begin{equation}\label{est:L22}
\sum_{|I|\le K-1}\left\|\frac{\partial \Gamma^{I}\phi}{(t-r)^{\frac{1}{2}}}\right\|_{L_{x}^{2}}\lesssim 
\sum_{|I|\le K-1}\left\|\partial \Gamma^{I}\phi\right\|_{L_{x}^{2}}\lesssim C\epsilon.
	\end{equation}
	We see that~\eqref{est:Keyphi} follows from~\eqref{est:L21} and~\eqref{est:L22}.
	
	\smallskip
	\textbf{Step 4.} Conclusion. Combining~\eqref{est:L2phi0} and~\eqref{est:Keyphi} with a standard induction argument, we have proved~\eqref{est:L2phi}.
	\end{proof}
\subsection{End of the proof of Theorem~\ref{thm:massless}}\label{SS:EndThm1}

We are in a position to complete the proof of Theorem~\ref{thm:massless} and start with the proof of Proposition~\ref{pro:main1}.

\begin{proof}[Proof of Proposition~\ref{pro:main1}]
	For any initial data $(\psi_{0},\vec{\phi}_{0})$ satisfying the smallness condition~\eqref{est:small}, we consider the corresponding solution $(\psi,\phi)$ of~\eqref{equ:3DDKG}. From the smallness condition~\eqref{est:small} and ${\rm{supp}}(\psi,\vec{\phi}_{0})\subseteq \mathcal{K}$, we know that 
	\begin{equation}\label{est:init}
	\mathcal{E}_{D}^{\delta}(2,\psi)^{\frac{1}{2}}+\mathcal{E}(2,\phi)^{\frac{1}{2}}+\mathcal{F}(2,\phi)^{\frac{1}{2}}\lesssim \epsilon.
	\end{equation}
	
	In what follows, we prove Proposition~\ref{pro:main1} by improving all the estimates of $(\psi,\phi)$ in the bootstrap assumption~\eqref{est:Boot1}.
	
	\smallskip
	\textbf{Step 1.} Closing the pointwise estimate of $\psi$. First, from~\eqref{est:Ga} and~\eqref{est:Psi}, we have 
	\begin{equation*}
	\sum_{|I|\le N-4}\left|\partial \Gamma^{I}\Psi\right|\lesssim (t-r)^{-1}\sum_{|I|\le N-3}\left| \Gamma^{I}\Psi\right|\lesssim \left(\epsilon+ C^{2}\epsilon^{2}\right)(t+r)^{-1+\delta}(t-r)^{-\frac{1}{2}-6\delta}.
	\end{equation*}
	Based on the above inequality,~$-i\gamma^{\mu}\partial_{\mu}\Psi=\psi$,~\eqref{est:pGammaGammahat} and ~\eqref{equ:comm}, we obtain
	\begin{equation*}
	\sum_{|I|\le N-4}\big|\widehat{\Gamma}^{I}\psi\big|\lesssim 	\sum_{|I|\le N-4}\left|\partial \Gamma^{I}\Psi\right|\lesssim \left(\epsilon+C^{2}\epsilon^{2}\right)(t+r)^{-1+\delta}(t-r)^{-\frac{1}{2}-6\delta}.
	\end{equation*}
	This inequality strictly improves the pointwise estimates of $\psi$ in the bootstrap assumption~\eqref{est:Boot1} for $C$ large enough and $\epsilon$ small enough.
	
	\smallskip
	\textbf{Step 2.}  Closing the pointwise estimate of $\phi$. Recall that, 	we decompose the wave component $\phi$ as,
	\begin{equation*}
	\phi(t,x)=\phi_{\emph{{free}}}(t,x)+\phi_{inho}(t,x)
	\end{equation*}
	where $\phi_{free}$ and $\phi_{inho}$ are the solutions for the following 3D linear homogeneous
	or inhomogeneous wave equations,
	\begin{equation*}
	\left\{\begin{aligned}
	-\Box \phi_{free}&=0, \quad \quad \ \ \mbox{with}\ \ (\phi_{free},\partial_{t}\phi_{free})_{|t=2}=(\phi_{0},\phi_{1}),\\
	-\Box \phi_{inho}&=\psi^{*}\gamma^{0}\psi,\quad  \mbox{with}\ \ (\phi_{inho},\partial_{t}\phi_{inho})_{|t=2}=(0,0).
	\end{aligned}\right.
	\end{equation*}
	
	First, from Lemma~\ref{le:waveenergy} and~\eqref{est:init}, we have 
	\begin{equation*}
	\sum_{|I|\le N}\left\|\Gamma^{I}\phi_{free}\right\|_{L_{x}^{2}}\lesssim \sum_{|I|\le N-1}\left(\mathcal{E}(t,\Gamma^{I}\phi_{free})^{\frac{1}{2}}+\mathcal{F}(t,\Gamma^{I}\phi_{free})^{\frac{1}{2}}\right)\lesssim \epsilon.
	\end{equation*}
	It follows from~\eqref{est:Klainer} that
	\begin{equation}\label{est:phiend}
	\sum_{|I|\le N-2}\left|\Gamma^{I}\phi_{free}\right|\lesssim \epsilon(t+r)^{-1}(t-r)^{-\frac{1}{2}}.
	\end{equation}
	
	Second, using~\eqref{est:Hidden},~\eqref{est:Gaga0},~\eqref{est:Boot1} and~\eqref{est:Firstpsi-}, we have 
	\begin{equation*}
	\begin{aligned}
	\sum_{|I|\le N-4}\left|\Gamma^{I}\left(\psi^{*}\gamma^{0}\psi\right)\right|
	&\lesssim \sum_{\substack{|I_{1}|\le N-4\\ |I_{2}|\le N-4}}\big|\widehat{\Gamma}^{I_{1}}\psi\big|\big|\big[\widehat{\Gamma}^{I_{2}}\psi\big]_{-}\big|\lesssim C^{2}\epsilon^{2}(t+r)^{-3+2\delta}(t-r)^{-10\delta}.
	\end{aligned}
	\end{equation*}
	On the other hand, from $N\ge 7$,~\eqref{est:Hidden},~\eqref{est:Gaga0},~\eqref{est:Boot1},~\eqref{est:Firstpsi},~\eqref{est:Firstpsi-} and~\eqref{est:Secondpsi-}, 
	\begin{equation*}
	\begin{aligned}
	\sum_{|I|\le N-3}\left|\Gamma^{I}\left(\psi^{*}\gamma^{0}\psi\right)\right|
	&\lesssim \sum_{\substack{|I_{1}|\le N-4\\ |I_{2}|\le N-3}}\big|\widehat{\Gamma}^{I_{1}}\psi\big|\big|\big[\widehat{\Gamma}^{I_{2}}\psi\big]_{-}\big|\\
	&+ \sum_{\substack{|I_{1}|\le N-3\\ |I_{2}|\le N-4}}\big|\widehat{\Gamma}^{I_{1}}\psi\big|\big|\big[\widehat{\Gamma}^{I_{2}}\psi\big]_{-}\big|
	\lesssim C^{2}\epsilon^{2}(t+r)^{-3+2\delta}(t-r)^{-4\delta}.
	\end{aligned}
	\end{equation*}
	Therefore, using again (iii) of Proposition~\ref{pro:main1},~\eqref{est:small} and~\eqref{est:Klainer}, we obtain
	\begin{equation}\label{est:Endphi}
	\begin{aligned}
	\sum_{|I|\le N-4}\left|\Gamma^{I}\phi_{inho}\right|&\lesssim \epsilon(t+r)^{-1}(t-r)^{-\frac{1}{2}}+C^{2}\epsilon^{2} (t+r)^{-1}(t-r)^{-8\delta},\\
	\sum_{|I|\le N-3}\left|\Gamma^{I}\phi_{inho}\right|&\lesssim \epsilon(t+r)^{-1}(t-r)^{-\frac{1}{2}}+C^{2}\epsilon^{2} (t+r)^{-1}(t-r)^{-2\delta}.
	\end{aligned}
	\end{equation}
	
	Combining~\eqref{est:phiend} and~\eqref{est:Endphi}, we obtain
	\begin{equation*}
	\begin{aligned}
	\sum_{|I|\le N-4}\left|\Gamma^{I}\phi\right|&\lesssim \left(\epsilon+C^{2}\epsilon^{2}\right) (t+r)^{-1}(t-r)^{-8\delta},\\
	\sum_{|I|\le N-3}\left|\Gamma^{I}\phi\right|&\lesssim \left(\epsilon+C^{2}\epsilon^{2}\right) (t+r)^{-1}(t-r)^{-2\delta}.
	\end{aligned}
	\end{equation*}
	
	These two inequalities strictly improve the pointwise estimates of $\phi$ in~\eqref{est:Boot1} for $C$ large enough and $\epsilon$ small enough (depending on $C$).
	
	\smallskip
	\textbf{Step 3.} Bound on the energy norm of $\psi$. From~\eqref{est:GammaGammahat},~\eqref{est:Gamma-},~\eqref{est:Hidden},~\eqref{est:Gaga0}, $N\ge 7$ and an elementary computation, we have 
	\begin{equation*}
	\begin{aligned}
	\sum_{|I|\le N}\sum_{|J|\le N}\int_{2}^{t}\int_{\R^{3}}\big|\big(\widehat{\Gamma}^{I}\psi\big)^{*}\gamma^{0}\widehat{\Gamma}^{J}(\phi\psi)\big|\d x \d \tau\lesssim \mathcal{G}_{3}+\mathcal{G}_{4}+\mathcal{G}_{5}, 
	\end{aligned}
	\end{equation*}
	where
	\begin{equation*}
	\begin{aligned}
	\mathcal{G}_{3}&=\sum_{|I|\le N}\sum_{|I_{1}|\le N-3}\sum_{|I_{2}|\le N}\int_{2}^{t}\int_{\R^{3}}\big|\big[\widehat{\Gamma}^{I}\psi\big]_{-}\big|\big|\Gamma^{I_{1}}\phi\big|\big|\widehat{\Gamma}^{I_{2}}\psi\big|\d x \d \tau,\\
	\mathcal{G}_{4}&=\sum_{|I|\le N}\sum_{|I_{1}|\le N}\sum_{|I_{2}|\le N-4}\int_{2}^{t}\int_{\R^{3}}\big|\big[\widehat{\Gamma}^{I}\psi\big]_{-}\big|\big|\Gamma^{I_{1}}\phi\big|\big|\widehat{\Gamma}^{I_{2}}\psi\big|\d x \d \tau,\\
	\mathcal{G}_{5}&=\sum_{|I|\le N}\sum_{|I_{1}|\le N}\sum_{|I_{2}|\le N-4}\int_{2}^{t}\int_{\R^{3}}\big|\widehat{\Gamma}^{I}\psi\big|\big|\Gamma^{I_{1}}\phi\big|\big|\big[\widehat{\Gamma}^{I_{2}}\psi\big]_{-}\big|\d x \d \tau.
	\end{aligned}
	\end{equation*}
	
	First, from~\eqref{est:Boot1} and~\eqref{est:Ghostpsi}, we have 
	\begin{equation*}
	\begin{aligned}
	\mathcal{G}_{3}
	&\lesssim \sum_{|I|\le N}\sum_{\substack{|I_{1}|\le N-3\\ |I_{2}|\le N}}\int_{2}^{t}
	\left\|\frac{\big[\widehat{\Gamma}^{I}\psi\big]_{-}}{(\tau-r)^{\frac{1}{2}+\delta}}\right\|_{L_{x}^{2}}\big\|(\tau-r)^{\frac{1}{2}+\delta}\Gamma^{I_{1}}\phi\big\|_{L_{x}^{\infty}}\big\|\widehat{\Gamma}^{I_{2}}\psi\big\|_{L_{x}^{2}}\d \tau\\
	&\lesssim C^{2}\epsilon^{2}\sum_{|I|\le N}\int_{2}^{t}\tau^{-\frac{1}{2}-\delta}\left\|\frac{\big[\widehat{\Gamma}^{I}\psi\big]_{-}}{(\tau-r)^{\frac{1}{2}+\delta
	}}\right\|_{L_{x}^{2}}\d \tau \lesssim C^{3}\epsilon^{3}.
	\end{aligned}
	\end{equation*}
	
	Second, from~\eqref{est:Boot1},~\eqref{est:Ghostpsi} and~\eqref{est:L2phi}, we find
	\begin{equation*}
	\begin{aligned}
	\mathcal{G}_{4}
	&\lesssim \sum_{|I|\le N}\sum_{\substack{|I_{1}|\le N\\ |I_{2}|\le N-4}}\int_{2}^{t}
	\left\|\frac{\big[\widehat{\Gamma}^{I}\psi\big]_{-}}{(\tau-r)^{\frac{1}{2}+\delta}}\right\|_{L_{x}^{2}}\left\|\frac{{\Gamma}^{I_{1}}\phi}{(\tau-r)^{\frac{1}{2}}}\right\|_{L_{x}^{2}}
	\big\|(\tau-r)^{1+\delta}\widehat{\Gamma}^{I_{2}}\psi\big\|_{L_{x}^{\infty}}\d \tau\\
	&\lesssim C^{2}\epsilon^{2}\sum_{|I|\le N}\int_{2}^{t}\tau^{-\frac{1}{2}-\delta}\left\|\frac{\big[\widehat{\Gamma}^{I}\psi\big]_{-}}{(\tau-r)^{\frac{1}{2}+\delta}}\right\|_{L_{x}^{2}}\d \tau \lesssim C^{3}\epsilon^{3}.
	\end{aligned}
	\end{equation*}
	Then, using again~\eqref{est:Boot1},~\eqref{est:Ghostpsi},~\eqref{est:Firstpsi-} and~\eqref{est:L2phi},
	\begin{equation*}
	\begin{aligned}
	\mathcal{G}_{5}
	&\lesssim \sum_{|I|\le N}\sum_{\substack{|I_{1}|\le N\\ |I_{2}|\le N-4}}
	\int_{2}^{t}\big\|\widehat{\Gamma}^{I}\psi\big\|_{L_{x}^{2}}\left\|\frac{{\Gamma}^{I_{1}}\phi}{(\tau-r)^{\frac{1}{2}}}\right\|_{L_{x}^{2}}\big\|(\tau -r)^{\frac{1}{2}}\big[\widehat{\Gamma}^{I_{2}}\psi\big]_{-}\big\|_{L_{x}^{2}}\d \tau\\
	&\lesssim C^{3}\epsilon^{3} \int_{2}^{t}\tau^{-1-4\delta}\d \tau \lesssim C^{3}\epsilon^{3}.
	\end{aligned}
	\end{equation*}
	Combining the above inequalities, one can obtain
	\begin{equation*}
		\sum_{|I|\le N}\sum_{|J|\le N}\int_{2}^{t}\int_{\R^{3}}\big|\big(\widehat{\Gamma}^{I}\psi\big)^{*}\gamma^{0}\widehat{\Gamma}^{J}(\phi\psi)\big|\d x \d \tau\lesssim \mathcal{G}_{3}+\mathcal{G}_{4}+\mathcal{G}_{5}\lesssim C^{3}\epsilon^{3}.
	\end{equation*}
	
Therefore, from Lemma~\ref{le:DiracGhotst},~\eqref{equ:3DDKGuse} and~\eqref{est:init}, we obtain
	\begin{equation*}
	\begin{aligned}
	\sum_{|I|\le N}\mathcal{E}_{D}^{\delta}\big(t,\widehat{\Gamma}^{I}\psi\big)
	&\lesssim \sum_{|I|\le N}\sum_{|J|\le N}\int_{2}^{t}\int_{\R^{3}}\big|\big(\widehat{\Gamma}^{I}\psi\big)^{*}\gamma^{0}\widehat{\Gamma}^{J}(\phi\psi)\big|\d x \d \tau\\
	&+\sum_{|I|\le N}\mathcal{E}_{D}^{\delta}\big(2,\widehat{\Gamma}^{I}\psi\big)\lesssim \epsilon^{2}+C^{3}\epsilon^{3}.
	\end{aligned}
	\end{equation*}
	This inequality strictly improves the energy estimates of $\psi$ in the bootstrap assumption~\eqref{est:Boot1} for $C$ large enough and $\epsilon$ small enough.

	\smallskip
	\textbf{Step 4.} Bound on the energy norm of $\phi$. From $N\ge 7$,~\eqref{est:Hidden} and~\eqref{est:Gaga0}, we have 
	\begin{equation*}
	\sum_{|I|\le N}\int_{2}^{t}\left\|{\Gamma}^{I}\big(\psi^{*}\gamma^{0}\psi\big)\right\|_{L_{x}^{2}}\d \tau \lesssim \mathcal{G}_{6}+\mathcal{G}_{7},
	\end{equation*}
	where 
	\begin{equation*}
	\begin{aligned}
	\mathcal{G}_{6}&=\sum_{|I_{1}|\le N}\sum_{|I_{2}|\le N-2}\int_{2}^{t}\big\| \big|\big[\widehat{\Gamma}^{I_{1}}\psi\big]_{-}\big|\widehat{\Gamma}^{I_{2}}\psi\big|\big\|_{L_{x}^{2}}\d \tau,\\
	\mathcal{G}_{7}&=\sum_{|I_{1}|\le N-3}\sum_{|I_{2}|\le N}\int_{2}^{t}\big\| \big|\big[\widehat{\Gamma}^{I_{1}}\psi\big]_{-}\big|\widehat{\Gamma}^{I_{2}}\psi\big|\big\|_{L_{x}^{2}}\d \tau.
	\end{aligned}
	\end{equation*}
	
	From~\eqref{est:Boot1},~\eqref{est:Ghostpsi} and~\eqref{est:Firstpsi-}, we have 
	\begin{equation*}
	\begin{aligned}
	\mathcal{G}_{6}&\lesssim C\epsilon\sum_{|I_{1}|\le N}\int_{2}^{t}\tau^{-1+\delta}\left\|\frac{\big[\widehat{\Gamma}^{I_{1}}\psi\big]_{-}}{(\tau-r)^{\frac{1}{2}+\delta}}\right\|_{L_{x}^{2}}\d \tau\lesssim C^{2}\epsilon^{2},\\
	\mathcal{G}_{7}&\lesssim \sum_{|I_{1}|\le N-3}\sum_{|I_{2}|\le N}\int_{2}^{t}\big\| \big[\widehat{\Gamma}^{I_{1}}\psi\big]_{-}\big\|_{L_{x}^{\infty}}\big\|\widehat{\Gamma}^{I_{2}}\psi\big\|_{L_{x}^{2}}\d \tau\lesssim C^{2}\epsilon^{2}.
	\end{aligned}
	\end{equation*}
	Combining the above inequalities with~\eqref{equ:3DDKGuse},~\eqref{est:init} and Lemma~\ref{le:waveenergy}, we obtain
	\begin{equation*}
	\begin{aligned}
	\sum_{|I|\le N}\left\|\partial \Gamma^{I}\phi\right\|_{L_{x}^{2}}
	\lesssim \sum_{|I|\le N}\left(\mathcal{E}(2,\Gamma^{I}\phi)^{\frac{1}{2}}+\int_{2}^{t}\left\|{\Gamma}^{I}\big(\psi^{*}\gamma^{0}\psi\big)\right\|_{L_{x}^{2}}\d \tau\right)\lesssim \epsilon+C^{2}\epsilon^{2}.
	\end{aligned}
	\end{equation*}
	This inequality strictly improves the energy estimates of $\phi$ in the bootstrap assumption~\eqref{est:Boot1} for $C$ large enough and $\epsilon$ small enough.

	\smallskip
	At this point, we have strictly improved all the pointwise and energy estimates of $(\psi,\phi)$ in the bootstrap assumption~\eqref{est:Boot1}.  In conclusion, for all initial data $(\psi_{0},\vec{\phi}_{0})$ satisfying~\eqref{est:small}, we show that $T_{*}(\psi_{0},\vec{\phi}_{0})=\infty$ and thus the proof of Proposition~\ref{pro:main1} is complete.
	\end{proof}

Now, we complete the proof of Theorem~\ref{thm:massless} by using Proposition~\ref{pro:main1}.

\begin{proof}[End of the proof of Theorem~\ref{thm:massless}]
	Notice that, from~\eqref{est:Boot1}, the global
	existence of $(\psi,\phi)$, the pointwise decay~\eqref{est:theorem1point} and the energy bound~\eqref{est:theorem1} are consequences of Proposition~\ref{pro:main1}.
	Our task now is to prove that the wave component $\phi$ scatters to a free solution in the
	energy space $\mathcal{X}_{N}=\big(\dot{H}^{N+1}\cap \dot{H}^{1}\big)\times \big(\dot{H}^{N}\cap L^{2}\big)$. 
	
	\smallskip
	In order to discuss the scattering theory for the wave component $\phi$, we first recall a few facts about the 3D linear wave equation. Consider the Cauchy problem for the 3D linear wave equation, 
	\begin{equation}\label{equ:wave}
	-\Box u=G\quad \mbox{with}\quad (u,\pt u)_{|t=2}=(u_{0},u_{1}).
	\end{equation}
	The Fourier analysis gives us a solution of~\eqref{equ:wave},
	\begin{equation}\label{equ:Duhamwave}
	\begin{pmatrix}
	u\\ \pt u
	\end{pmatrix}(t)=\mathcal{S}_{W}(t-2)
	\left(\begin{aligned}
	u_{0}\\ u_{1}
	\end{aligned}\right)+\int_{2}^{t}\mathcal{S}_{W}(t-\tau)
	\left(\begin{aligned}
	&0\\ G&(\tau)
	\end{aligned}\right)\d \tau,
	\end{equation}
	where $\mathcal{S}_{W}(t)$ can be written as 
	\begin{equation*}
	\mathcal{S}_{W}(t)=
	\begin{pmatrix}
	\cos \left(t\sqrt{-\Delta}\right)& (\sqrt{-\Delta})^{-1}\sin  \left(t\sqrt{-\Delta}\right)\\
	-\sqrt{-\Delta}\sin  \left(t\sqrt{-\Delta}\right)& \cos \left(t\sqrt{-\Delta}\right)
	\end{pmatrix},\quad \forall t\in (-\infty,\infty).
	\end{equation*}
	Here, the operators mean that 
	\begin{equation*}
	\begin{aligned}
	\left((\sqrt{-\Delta})^{-1}f\right)(x)&=\mathcal{F}^{-1}\left(|\xi|^{-1}\hat{f}(\xi)\right)(x),\\
	\left({\cos(t\sqrt{-\Delta})f}\right)(x)&=\mathcal{F}^{-1}\left(\cos (t|\xi|)\hat{f}(\xi)\right)(x),
	\end{aligned}
	\end{equation*}
	and other operators are similarly defined. Recall that, for $K\in \mathbb{N}^{+}$, the 3D wave flow $\mathcal{S}_{W}(t)$ are isometric groups on $\dot{H}^{K}\times\dot{H}^{K-1}$. More precisely, for $K\in \mathbb{N}^{+}$, we have the following identities,
	\begin{equation}\label{est:Sw}
	\left\|\mathcal{S}_{W}(t)\right\|_{\dot{H}^{K}\times \dot{H}^{K-1}}=1,\quad \forall t\in (-\infty,\infty).
	\end{equation}
	Recall also that, we also have the following basic properties of $\mathcal{S}_{W}$,
	\begin{equation}\label{equ:Sts}
	\mathcal{S}_{W}(t+\tau)=\mathcal{S}_{W}(t)\circ\mathcal{S}_{W}(\tau),\quad \forall(t,\tau)\in (-\infty,\infty)^{2}.
	\end{equation}
	
	\smallskip
	Now, we show that the wave component $\phi$ scatters to a free solution. Actually, from the Step 4 of the proof of Proposition~\ref{pro:main1} and $T_{*}(\psi_{0},\vec{\phi}_{0})=\infty$, we know that 
	\begin{equation}\label{est:psiHN}
	\sum_{|I|\le N}\int_{2}^{\infty}\left\|\Gamma^{I}(\psi^{*}\gamma^{0}\psi)\right\|_{L_{x}^{2}}\d \tau\lesssim C\epsilon\Rightarrow
\int_{2}^{\infty}\left\|\psi^{*}\gamma^{0}\psi\right\|_{H^{N}}\d \tau\lesssim C\epsilon.
	\end{equation}
	Based on the above inequality and~\eqref{est:Sw}, we have 
	\begin{equation*}
	\int_{2}^{\infty}\left\|\mathcal{S}_{W}(2-\tau)	\left(\begin{aligned}
	&0\\ \psi^{*}&\gamma^{0}\psi
	\end{aligned}\right)\right\|_{\mathcal{X}_{N}}\d \tau \lesssim \int_{2}^{\infty}\left\|\psi^{*}\gamma^{0}\psi\right\|_{H^{N}}\d \tau\lesssim C\epsilon,
	\end{equation*}
	which implies
	\begin{equation*}
	\begin{aligned}
	\begin{pmatrix}
	\phi_{0h}\\ \phi_{1h}
	\end{pmatrix}:&=\int_{2}^{\infty}\mathcal{S}_{W}(2-\tau)\begin{pmatrix}
	0\\
	\psi^{*}\gamma^{0}\psi
	\end{pmatrix}\d \tau \in \mathcal{X}_{N}.
	\end{aligned}
	\end{equation*}
	Therefore,
		\begin{equation*}
	\begin{pmatrix}
	\phi_{0\ell}\\ \phi_{1\ell}
	\end{pmatrix}:=\begin{pmatrix}
	\phi_{0}\\ \phi_{1}
	\end{pmatrix}
	+\begin{pmatrix}
	\phi_{0h}\\ \phi_{1h}
	\end{pmatrix}\in \mathcal{X}_{N}\Rightarrow \begin{pmatrix}
	\phi_{\ell}\\ \pt \phi_{\ell}
	\end{pmatrix}(t):=\mathcal{S}_{W}(t-2)	\begin{pmatrix}
	\phi_{0\ell}\\ \phi_{1\ell}
	\end{pmatrix}\in \mathcal{X}_{N}.
	\end{equation*}
	Last, from~\eqref{est:Sw},~\eqref{equ:Sts} and~\eqref{est:psiHN}, we obtain
	\begin{equation*}
	\lim_{t\to \infty}\left\|(\phi,\pt \phi)-(\phi_{\ell},\pt \phi_{\ell})\right\|_{\mathcal{X}_{N}}\le \lim_{t\to \infty}\int_{t}^{\infty}\left\|\psi^{*}\gamma^{0}\psi\right\|_{H^{N}}\d \tau=0,
	\end{equation*}
	which means that the wave component $\phi$ scatters to a free solution in the energy space $\mathcal{X}_{N}$. The proof of Theorem~\ref{thm:massless} is complete.
	\end{proof}

\section{Proof of Theorem~\ref{thm:massive}}\label{S:thmmassive}

In this section, we prove the existence of global-in-time solution $(\psi,\phi)$ of~\eqref{equ:3DDKG2} satisfying~\eqref{est:theorem2point}-\eqref{est:theorem2sca} in Theorem~\ref{thm:massive} via a bootstrap argument. For simplicity of notation, we denote $\Upsilon=i\gamma^{0}\gamma^{5}=-\gamma^{1}\gamma^{2}\gamma^{3}$ in the rest of the section.

\subsection{Bootstrap assumption}
Let $N\in \mathbb{N}^{+}$ with $N\ge 13$ and $0<\delta\ll 1$. The proof of Theorem~\ref{thm:massive} relies on a bootstrap argument of high-order energy and pointwise decay of solution $(\psi,\phi)$. More precisely, we consider the following bootstrap assumption: for $C_{0}\gg 1$ and $0<\varepsilon\ll C_{0}^{-1}$ to be chosen later,

\begin{equation}\label{est:Boot21}
\left\{
\begin{aligned}\sum_{|I|\le N}\left(\mathcal{E}_{D}^{\delta}(t,\widehat{Z}^{I}\psi)^{\frac{1}{2}}+\mathcal{E}(t,Z^{I}\phi)^{\frac{1}{2}}\right)&\le C_{0}\varepsilon,\\
\sum_{|I|\le N-6}\langle t+r\rangle^{\frac{3}{2}}\left|Z^{I}\psi\right|+\sum_{|I|\le N-2}\left\|\frac{\langle t+r\rangle }{\langle t-r\rangle}Z^{I}\psi\right\|_{L_{x}^{2}}&\le C_{0}\varepsilon,\\
\sum_{|I|\le N-7}\langle t+r\rangle \langle t-r\rangle^{\frac{1}{2}}\left|Z^{I}\phi\right|+\sum_{|I|\le N}\left\|Z^{I}\phi\right\|_{L_{x}^{2}}&\le C_{0}\varepsilon.
\end{aligned}
\right.
\end{equation}

For all initial data $(\psi_{0},\vec{\phi}_{0})$ satisfying~\eqref{est:small2}, we set 
\begin{equation}\label{def:T*2}
T^{*}:=T^{*}(\psi_{0},\vec{\phi}_{0})=\sup \left\{t\in [2,\infty):(\psi,\phi) \ \mbox{satisfy~\eqref{est:Boot21}}\ \mbox{on}\ [2,t)\right\}>2.
\end{equation}

The following proposition is the main part of the proof of Theorem~\ref{thm:massive}.
\begin{proposition}\label{pro:main2}
	For any initial data $(\psi_{0},\vec{\phi}_{0})$ satisfying the smallness condition~\eqref{est:small2} in Theorem~\ref{thm:massive}, we have $T^{*}(\psi_{0},\vec{\phi}_{0})=\infty$.
\end{proposition}

The rest of the section is organized as follows. First, in \S\ref{SS:Keynonlinear}, we introduce several key spacetime norm estimates for the nonlinear terms. Later, in \S\ref{SS:EndThm2}, we prove Proposition~\ref{pro:main2} and then finish the proof of Theorem~\ref{thm:massive} from Proposition~\ref{pro:main2}. In what follows, the implied constants in $\lesssim$ do not depend on the
constants $C_{0}$ and $\varepsilon$ appearing in the bootstrap assumption~\eqref{est:Boot21}.

\subsection{Key estimate for the nonlinear terms}\label{SS:Keynonlinear}
In this subsection, we establish the estimates on the nonlinear terms that will be needed in the proof of Proposition~\ref{pro:main2}. 

From~\eqref{equ:3DDKG2} and~\eqref{equ:comm}, for $I\in \mathbb{N}^{10}$ with $|I|\le N$, we have  
\begin{equation}\label{equ:3DDKG2vec}
-i\gamma^{\mu}\partial_{\mu}\widehat{Z}^{I}\psi+\widehat{Z}^{I}\psi=\widehat{Z}^{I}\left(i\phi\gamma^{5}\psi\right)\quad \mbox{and}\quad 
-\Box Z^{I}\phi=Z^{I}\big(i\psi^{*}\gamma^{0}\gamma^{5}\psi\big),
\end{equation}
which will be frequently used in this section.

Note that, from~\eqref{est:parGamma},~\eqref{est:SobolevGlobal},~\eqref{est:Boot21} and the Cauchy-Schwarz inequality, we have 
\begin{equation}\label{est:usefulpsiphi}
\begin{aligned}
\sum_{|I|\le N}\int_{2}^{t}\tau^{-\frac{1}{2}-\delta}\left\|\frac{\big[\widehat{Z}^{I}\psi\big]_{-}}{\langle \tau -r\rangle^{\frac{1}{2}+\delta}}\right\|_{L_{x}^{2}}\d \tau&\lesssim C_{0}\varepsilon,\\
(t+r)^{\frac{3}{4}}\sum_{|I|\le N-3}\left( \big|\widehat{Z}^{I}\psi\big|+\left|Z^{I}\phi\right|+\left|\partial Z^{I}\phi\right|\right)&\lesssim C_{0}\varepsilon.
\end{aligned}
\end{equation}

Note also that, from the definition of $\gamma^{5}=i\gamma^{0}\gamma^{1}\gamma^{2}\gamma^{3}$ and~\eqref{equ:gamma}, we have
	\begin{equation}\label{equ:gamma5}
	(\gamma^{5})^{*}=\gamma^{5},\quad (\gamma^{5})^{2}={\rm{Id}},\quad \left\{\gamma^{5},\gamma^{\alpha}\right\}=0,\quad [\gamma^{5}\varphi]_{-}=\gamma^{5}[\varphi]_{-},
	\end{equation}
	where $\alpha\in \left\{0,1,2,3\right\}$ and $\varphi=\varphi(t,x):\R^{1+3}\to \C^{4}$.

\smallskip
Now we establish the $L_{t}^{1}L_{x}^{2}$ estimates for the nonlinear terms that will be needed in the energy estimates of Dirac spinor $\psi$.

\begin{lemma}
	For all $t\in [2,T^{*}(\psi_{0},\vec{\phi}_{0}))$, we have 
	\begin{equation}\label{est:Lt1Lx2psi}
	\sum_{|I|\le N}\int_{2}^{t}\int_{\R^{3}}\big|\big(\widehat{Z}^{I}\psi\big)\gamma^{0}\widehat{Z}^{I}\big(i\phi \gamma^{5}\psi\big)\big|\d x\d \tau\lesssim C_{0}^{3}\varepsilon^{3}.
	\end{equation}
\end{lemma}
\begin{proof}
	Note that, from~\eqref{est:pGammaGammahat},~\eqref{est:Z-},~\eqref{est:Hidden},~\eqref{equ:gamma5} and $N\ge 13$, we have 
	\begin{equation*}
	\sum_{|I|\le N}\int_{2}^{t}\int_{\R^{3}}\big|\big(\widehat{Z}^{I}\psi\big)\gamma^{0}\widehat{Z}^{I}\big(i\phi \gamma^{5}\psi\big)\big|\d x\d \tau\lesssim \mathcal{J}_{1}+\mathcal{J}_{2},
	\end{equation*}
	where
	\begin{equation*}
	\begin{aligned}
	\mathcal{J}_{1}&=\sum_{|I|\le N}\sum_{|I_{1}|\le N}\sum_{|I_{2}|\le N-6}\int_{2}^{t}\int_{\R^{3}}\big|\widehat{Z}^{I}\psi\big|\left|Z^{I_{1}}\phi\right|\big|{Z}^{I_{2}}\psi\big|\d x \d \tau,\\
	\mathcal{J}_{2}&=\sum_{|I|\le N}\sum_{|I_{1}|\le N-7}\sum_{|I_{2}|\le N}\int_{2}^{t}\int_{\R^{3}}\big|\widehat{Z}^{I}\psi\big|\left|Z^{I_{1}}\phi\right|\big|\big[\widehat{Z}^{I_{2}}\psi\big]_{-}\big|\d x \d \tau.
	\end{aligned}
	\end{equation*}
	From~\eqref{est:Boot21} and~\eqref{est:usefulpsiphi}, we see that 
	\begin{equation*}
	\begin{aligned}
		\mathcal{J}_{1}&\lesssim \sum_{|I|\le N}\sum_{\substack{|I_{1}|\le N\\|I_{2}|\le N-6}}\int_{2}^{t}
	\big\|\widehat{Z}^{I}\psi\big\|_{L_{x}^{2}}\left\|Z^{I_{1}}\phi\right\|_{L_{x}^{2}}\big\|{Z}^{I_{2}}\psi\big\|_{L_{x}^{\infty}} \d \tau\lesssim C_{0}^{3}\varepsilon^{3},\\
		\mathcal{J}_{2}
		&\lesssim\sum_{|I|\le N}\sum_{\substack{|I_{1}|\le N-7\\ |I_{2}|\le N}}\int_{2}^{t}\big\|\widehat{Z}^{I}\psi\big\|_{L_{x}^{2}}\left\|\langle \tau-r\rangle^{\frac{1}{2}+\delta}Z^{I_{1}}\phi\right\|_{L_{x}^{\infty}}\left\|\frac{\big[\widehat{Z}^{I_{2}}\psi\big]_{-}}{\langle \tau-r\rangle^{\frac{1}{2}+\delta}}\right\|_{L_{x}^{2}}\d \tau\\
		&\lesssim C_{0}^{2}\varepsilon^{2}\sum_{|I|\le N}\int_{2}^{t}\tau^{-1+\delta}\left\|\frac{\big[\widehat{Z}^{I}\psi\big]_{-}}{\langle \tau-r\rangle^{\frac{1}{2}+\delta}}\right\|_{L_{x}^{2}}\d \tau\lesssim C_{0}^{3}\varepsilon^{3}.
	\end{aligned}
	\end{equation*}
	
	Combining the above two inequalities, we obtain~\eqref{est:Lt1Lx2psi}.
	\end{proof}

Second, we establish the weighted $L_{x}^{2}$ estimate for the nonlinear terms which will be needed in the $L^{\infty}$ estimates of Dirac spinor $\psi$.

\begin{lemma}
	For all $t\in [2,T^{*}(\psi_{0},\vec{\phi}_{0}))$, we have 
	\begin{equation}\label{est:Lx2psi}
	\begin{aligned}
	\sum_{|I|\le N-2} \left\|(t+r)Z^{I}(\psi \partial \phi)\right\|_{L^{2}_{x}}
	&\lesssim C_{0}^{2}\varepsilon^{2}t^{-\frac{1}{2}},\\
		\sum_{|I|\le N-2} \left\|(t+r)Z^{I}(\phi^{2}\psi )\right\|_{L^{2}_{x}}
	&\lesssim C_{0}^{3}\varepsilon^{3}t^{-1}.
	\end{aligned}
	\end{equation}
\end{lemma}

\begin{proof}
	Note that, from~\eqref{est:GammaGammahat}, the Leibniz rule and $N\ge 13$, we have 
		\begin{equation*}
	\begin{aligned}
	\sum_{|I|\le N-2} \left\|(t+r)Z^{I}(\psi \partial \phi)\right\|_{L^{2}_{x}}
	&\lesssim \mathcal{J}_{3}+\mathcal{J}_{4},\\
	\sum_{|I|\le N-2} \left\|(t+r)Z^{I}(\phi^{2}\psi )\right\|_{L^{2}_{x}}
	&\lesssim \mathcal{J}_{5}+\mathcal{J}_{6},
	\end{aligned}
	\end{equation*}
	where
	\begin{equation*}
	\begin{aligned}
	\mathcal{J}_{3}&=\sum_{|I_{1}|\le N-6}\sum_{|I_{2}|\le N-1}\left\|(t+r)\left(Z^{I_{1}}\psi\right)\left(Z^{I_{2}}\phi\right)\right\|_{L_{x}^{2}},\\
	\mathcal{J}_{4}&=\sum_{|I_{1}|\le N-2}\sum_{|I_{2}|\le N-7}\left\|(t+r)\left(Z^{I_{1}}\psi\right)\left(Z^{I_{2}}\phi\right)\right\|_{L_{x}^{2}},\\
	\mathcal{J}_{5}&=\sum_{|I_{1}|\le N-2}\sum_{|I_{2}|\le N-7}\sum_{|I_{3}|\le N-7}\left\|(t+r)\left(Z^{I_{1}}\psi\right)\left(Z^{I_{2}}\phi\right)\left(Z^{I_{3}}\phi\right)\right\|_{L_{x}^{2}},\\
	\mathcal{J}_{6}&=\sum_{|I_{1}|\le N-6}\sum_{|I_{2}|\le N-7}\sum_{|I_{3}|\le N-2}\left\|(t+r)\left(Z^{I_{1}}\psi\right)\left(Z^{I_{2}}\phi\right)\left(Z^{I_{3}}\phi\right)\right\|_{L_{x}^{2}}.
	\end{aligned}
	\end{equation*}
	Using~\eqref{est:Boot21}, we have
	\begin{equation*}
	\begin{aligned}
	\mathcal{J}_{3}&\lesssim \sum_{\substack{|I_{1}|\le N-6\\|I_{2}|\le N-1}}\left\|(t+r)Z^{I_{1}}\psi\right\|_{L_{x}^{\infty}}\left\|Z^{I_{2}}\phi\right\|_{L_{x}^{2}}\lesssim C_{0}^{2}\varepsilon^{2}t^{-\frac{1}{2}},\\
	\mathcal{J}_{4}&\lesssim \sum_{\substack{|I_{1}|\le N-2\\|I_{2}|\le N-7}}\left\|\langle t-r\rangle Z^{I_{2}}\phi\right\|_{L_{x}^{\infty}}\left\|\frac{\langle t+r\rangle}{\langle t-r \rangle}Z^{I_{1}}\psi\right\|_{L_{x}^{2}}\lesssim C_{0}^{2}\varepsilon^{2}t^{-\frac{1}{2}},\quad \quad \quad \quad \quad 
	\end{aligned}
	\end{equation*}
		\begin{equation*}
	\begin{aligned}
	\mathcal{J}_{5}&\lesssim \sum_{|I_{1}|\le N-2}\sum_{\substack{|I_{2}|\le N-7\\|I_{3}|\le N-7}}\left\|Z^{I_{1}}\psi\right\|_{L_{x}^{2}}\left\|(t+r)\left(Z^{I_{2}}\phi\right)\left(Z^{I_{3}}\phi\right)\right\|_{L_{x}^{\infty}}\lesssim C^{3}_{0}\varepsilon^{3}t^{-1},\\
	\mathcal{J}_{6}&\lesssim \sum_{|I_{1}|\le N-6}\sum_{\substack{|I_{2}|\le N-7\\|I_{3}|\le N-2}}\left\|(t+r)\left(Z^{I_{1}}\psi\right)\left(Z^{I_{2}}\phi\right)\right\|_{L_{x}^{\infty}}\left\|Z^{I_{3}}\phi\right\|_{L_{x}^{2}}\lesssim C_{0}^{3}\varepsilon^{3}t^{-\frac{3}{2}}.
	\end{aligned}
	\end{equation*}
	Combining the above inequalities, we obtain~\eqref{est:Lx2psi}.
	\end{proof}

Third, following~\cite{Bache}, we deduce the hidden structure within the 3D Dirac--Klein-Gordon system~\eqref{equ:3DDKG2} for the wave component $\phi$. We set 
\begin{equation}\label{equ:N1N2}
\begin{aligned}
\mathcal{N}_{1}(\psi,\phi)&=-\phi\psi^{*}\gamma^{0}\left(i\gamma^{\mu}\partial_{\mu}\psi+\psi\right),\quad \mathcal{N}_{2}(\psi,\phi)=-\phi^{2}\psi^{*}\Upsilon\psi,\\
\mathcal{N}_{3}(\psi,\psi^{*})&=Q_{0a}\left(\psi,i\gamma^{5}\gamma^{a}\psi\right)+\sum_{b<a}Q^{ab}\left(\gamma^{0}\gamma^{a}\psi,i\gamma^{5}\gamma^{b}\psi\right).
\end{aligned}
\end{equation}

\begin{lemma}[\cite{Bache}]\label{le:Hidden}
	Let $\tilde{\phi}=\phi+\frac{1}{4}\psi^{*}\Upsilon\psi$, then we have
	\begin{equation}\label{equ:Hiddenphi}
	-\Box \tilde{\phi}=\frac{1}{2}\mathcal{N}_{1}(\psi,\phi)+\frac{1}{2}\mathcal{N}_{2}(\psi,\phi)+\frac{1}{2}\mathcal{N}_{3}(\psi,\psi
	^{*}).
	\end{equation}
\end{lemma}

\begin{proof}
	We divide the proof into three steps as follows.
	
	\smallskip
	\textbf{Step 1.} Expansion of $\psi^{*}\Upsilon\psi$. We claim that 
	\begin{equation}\label{equ:psiHpsi}
	\psi^{*}\Upsilon\psi=Q_{0}\left(\psi,\Upsilon\psi\right)+\mathcal{N}_{1}(\psi,\phi)+\mathcal{N}_{3}(\psi,\psi^{*}).
	\end{equation}
	
	Indeed, from~\eqref{equ:3DDKG2} and~\eqref{equ:gamma5}, we see that 
	\begin{equation*}
	\psi=i\gamma^{\mu}\partial_{\mu}\psi+i\phi\gamma^{5}\psi\quad \mbox{and}\quad 
	\psi^{*}=-i(\gamma^{\mu}\partial_{\mu}\psi)^{*}-i\phi\psi^{*}\gamma^{5}.
	\end{equation*}
	It follows from~\eqref{equ:gamma5} that
	\begin{equation*}
	\begin{aligned}
	\psi^{*}\Upsilon\psi
	&=\psi^{*}\Upsilon\left(i\gamma^{\mu}\partial_{\mu}\psi+i\phi\gamma^{5}\psi\right)\\
	&=\psi^{*}\Upsilon\left(i\gamma^{\mu}\partial_{\mu}\psi\right)-\phi\psi^{*}\gamma^{0}\psi=\left(\gamma^{\mu}\partial_{\mu}\psi\right)^{*}\Upsilon\left(\gamma^{\nu}\partial_{\nu}\psi\right)+\mathcal{N}_{1}\left(\psi,\phi\right).
	\end{aligned}
		\end{equation*}
	Using again~\eqref{equ:gamma5}, we find
	\begin{equation*}
	\left(\gamma^{\mu}\partial_{\mu}\psi\right)^{*}\Upsilon\left(\gamma^{\nu}\partial_{\nu}\psi\right)
	=\left[\left(\partial_{0}\psi\right)^{*}\gamma^{0}-\left(\partial_{a}\psi\right)^{*}\gamma^{a}\right]\Upsilon\left(\gamma^{0}\partial_{0}\psi+\gamma^{b}\partial_{b}\psi\right)=\mathcal{J}_{7}+\mathcal{J}_{8},
	\end{equation*}
	where (below the Einstein summation convention is not adopted when the indices are all at upper place)
	\begin{equation*}
\begin{aligned}
\mathcal{J}_{7}&=\partial_{0}\psi^{*}\gamma^{0}\Upsilon\gamma^{0}\partial_{0}\psi-\sum_{a=1}^{3}\partial^{a}\psi^{*}\gamma^{a}\Upsilon\gamma^{a}\partial^{a}\psi,\\
\mathcal{J}_{8}&=\partial_{0}\psi^{*}\gamma^{0}\Upsilon\gamma^{a}\partial_{a}\psi-\partial_{a}\psi^{*}\gamma^{a}\Upsilon\gamma^{0}\partial_{0}\psi-\sum_{a\ne b}\partial^{a}\psi^{*}\gamma^{a}\Upsilon\gamma^{b}\partial^{b}\psi.
	\end{aligned}
	\end{equation*}
	Moreover, from the definition of $\Upsilon$,~\eqref{equ:gamma} and~\eqref{equ:gamma5}, we find
	\begin{equation*}
	\mathcal{J}_{7}=-\partial_{0}\psi^{*}\partial_{0}\left(\Upsilon\psi\right)+\partial^{a}\psi^{*}\partial_{a}\left(\Upsilon\psi\right)=Q_{0}(\psi,\Upsilon\psi),
	\end{equation*}
	\begin{equation*}
	\begin{aligned}
	\mathcal{J}_{8}=
	&\sum_{a=1}^{3}\left(\partial_{0}\psi^{*}\partial^{a}\left(i\gamma^{5}\gamma^{a}\psi\right)-\partial^{a}\psi^{*}\partial_{0}\left(i\gamma^{5}\gamma^{a}\psi\right)\right)\\
	+&\sum_{b<a}\left(\partial^{a}\psi^{*}\gamma^{0}\gamma^{a}\partial^{b}\left(i\gamma^{5}\gamma^{b}\psi\right)-\partial^{b}\psi^{*}\gamma^{0}\gamma^{a}\partial^{a}\left(i\gamma^{5}\gamma^{b}\psi\right)\right)\\
	=&Q_{0a}\left(\psi,i\gamma^{5}\gamma^{a}\psi\right)+\sum_{b<a}Q^{ab}\left(\gamma^{0}\gamma^{a}\psi,i\gamma^{5}\gamma^{b}\psi\right).
	\end{aligned}
	\end{equation*}
Combining the above identities, we obtain~\eqref{equ:psiHpsi}.

	\smallskip
	\textbf{Step 2.} Expansion of $-\Box\left(\psi^{*}\Upsilon\psi\right)$. We claim that 
	\begin{equation}\label{equ:BoxpsiHpsi}
	-\Box \left(\psi^{*}\Upsilon\psi\right)=-2\psi^{*}\Upsilon\psi-2Q_{0}(\psi,\Upsilon\psi)+2\mathcal{N}_{2}(\psi,\phi).
	\end{equation}
	Indeed, by an elementary computation, we find
	\begin{equation*}
	-\Box \left(\psi^{*}\Upsilon\psi\right)=\left(-\Box \psi\right)^{*}\Upsilon\psi+\psi^{*}\Upsilon\left(-\Box \psi\right)-2Q_{0}\left(\psi,\Upsilon\psi\right).
	\end{equation*}
	Acting $i\gamma^{\mu} \partial_{\mu}$ on both sides of the first line of~\eqref{equ:3DDKG2} and then using $(i\gamma^{\mu}\partial_{\mu})^{2}=\Box$,
	\begin{equation*}
	\begin{aligned}
	-\Box\psi&=-\psi+i\phi\gamma^{5}\psi+i\gamma^{\mu}\partial_{\mu}\left(i\phi\gamma^{5}\psi\right)\\
	&=-\psi+i\gamma^{\mu}\partial_{\mu}\phi(i\gamma^{5}\psi)+i\phi\psi^{5}\left(\psi-i\gamma^{\mu}\partial_{\mu}\psi\right)\\
	&=-\psi-\phi^{2}\psi-(\partial_{\mu}\phi)\gamma^{\mu}\gamma^{5}\psi.
	\end{aligned}
	\end{equation*}
	On the other hand, from~\eqref{equ:gamma5}, we see that 
	\begin{equation*}
\Upsilon\gamma^{\mu}\gamma^{5}+\gamma^{5}\left(\gamma^{\mu}\right)^{*}\Upsilon=0,\quad \mbox{for all}\ \mu\in \left\{0,1,2,3\right\}.
	\end{equation*}
We see that~\eqref{equ:BoxpsiHpsi} follows from	the above three identities.
	
	\smallskip
	\textbf{Step 3.} Conclusion. Combining~\eqref{equ:psiHpsi}-\eqref{equ:BoxpsiHpsi} with the definition of $\tilde{\phi}$, we obtain
	\begin{equation*}
	\begin{aligned}
	-\Box \tilde{\phi}
	&=-\Box \phi-\frac{1}{4}\Box\left(\psi^{*}\Upsilon\psi\right)\\
	&=\frac{1}{2}\psi^{*}\Upsilon\psi-\frac{1}{2}Q_{0}\left(\psi,\Upsilon\psi\right)+\frac{1}{2}\mathcal{N}_{2}(\psi,\phi)\\
	&=\frac{1}{2}\mathcal{N}_{1}(\psi,\phi)+\frac{1}{2}\mathcal{N}_{2}(\psi,\phi)+\frac{1}{2}\mathcal{N}_{3}(\psi,\psi
	^{*}).
	\end{aligned}
	\end{equation*}
	The proof of~\eqref{equ:Hiddenphi} is complete.
	\end{proof}

Next, we establish the weighted $L_{t}^{1}L_{x}^{2}$ estimates for the nonlinear terms $\left(\mathcal{N}_{1},\mathcal{N}_{2},\mathcal{N}_{3}\right)$.

\begin{lemma}
	For all $t\in [2,T^{*}(\psi_{0},\vec{\phi}_{0}))$, we have
	\begin{equation}\label{est:N1N2}
	\begin{aligned}
	\sum_{|I|\le N-1}\int_{2}^{t}\left\|(\tau+r)Z^{I}\mathcal{N}_{1}(\psi,\phi)\right\|_{L_{x}^{2}}\d \tau&\lesssim C_{0}^{3}\varepsilon^{3},\\
	\sum_{|I|\le N-1}\int_{2}^{t}\left\|(\tau+r)Z^{I}\mathcal{N}_{2}(\psi,\phi)\right\|_{L_{x}^{2}}\d \tau&\lesssim C_{0}^{4}\varepsilon^{4},\\
	 \sum_{|I|\le N-1}\int_{2}^{t}\left\|(\tau+r)Z^{I}\mathcal{N}_{3}(\psi,\psi^{*})\right\|_{L_{x}^{2}}\d \tau&\lesssim C_{0}^{2}\varepsilon^{2}.
	 \end{aligned}
	\end{equation}
\end{lemma}

\begin{proof}
	Note that, from~\eqref{est:Qab}, the Leibniz rule and $N\ge 13$, we have 
	\begin{equation*}
	\begin{aligned}
		\sum_{|I|\le N-1} \int_{2}^{t}\left\|(\tau+r)Z^{I}\mathcal{N}_{3}(\psi,\psi^{*})\right\|_{L^{2}_{x}}\d \tau
	&\lesssim \mathcal{J}_{9},\\
	\sum_{|I|\le N-1} \int_{2}^{t}\left\|(\tau+r)Z^{I}\mathcal{N}_{1}(\psi,\phi)\right\|_{L^{2}_{x}}\d \tau
	&\lesssim \mathcal{J}_{10}+\mathcal{J}_{11},\\
	\sum_{|I|\le N-1} \int_{2}^{t}\left\|(\tau+r)Z^{I}\mathcal{N}_{2}(\psi,\phi)\right\|_{L^{2}_{x}}\d \tau
	&\lesssim \mathcal{J}_{12}+\mathcal{J}_{13},
	\end{aligned}
	\end{equation*}
	where
	\begin{equation*}
	\begin{aligned}
	\mathcal{J}_{9}&=\sum_{|I_{1}|\le N}\sum_{|I_{2}|\le N-6}\int_{2}^{t}\left\|\left(Z^{I_{1}}\psi\right)\left(Z^{I_{2}}\psi\right)\right\|_{L_{x}^{2}}\d \tau,\\
	\mathcal{J}_{10}&=\sum_{|I_{1}|\le N-1}\sum_{|I_{2}|\le N-6}\sum_{|I_{3}|\le N-6 }\int_{2}^{t}\left\|(\tau+r)\left(Z^{I_{1}}\phi\right)\left(Z^{I_{2}}\psi\right)\left(Z^{I_{3}}\psi\right)\right\|_{L_{x}^{2}}\d \tau,\\
	\mathcal{J}_{11}&=\sum_{|I_{1}|\le N-7}\sum_{|I_{2}|\le N-6}\sum_{|I_{3}|\le N }\int_{2}^{t}\left\|(\tau+r)\left(Z^{I_{1}}\phi\right)\left(Z^{I_{2}}\psi\right)\left(Z^{I_{3}}\psi\right)\right\|_{L_{x}^{2}}\d \tau,
		\end{aligned}
	\end{equation*}
		\begin{equation*}
	\begin{aligned}
	\mathcal{J}_{12}&=\sum_{\substack{|I_{1}|\le N-1\\ |I_{2}|\le N-7}}\sum_{\substack{|I_{3}|\le N-6\\|I_{4}|\le N-6}}\int_{2}^{t}\left\|(\tau+r)\left(Z^{I_{1}}\phi\right)\left(Z^{I_{2}}\phi\right)\left(Z^{I_{3}}\psi\right)\left(Z^{I_{4}}\psi\right)\right\|_{L_{x}^{2}}\d \tau,\\
	\mathcal{J}_{13}&=\sum_{\substack{|I_{1}|\le N-7\\ |I_{2}|\le N-7}}\sum_{\substack{|I_{3}|\le N-1\\|I_{4}|\le N-6}}\int_{2}^{t}\left\|(\tau+r)\left(Z^{I_{1}}\phi\right)\left(Z^{I_{2}}\phi\right)\left(Z^{I_{3}}\psi\right)\left(Z^{I_{4}}\psi\right)\right\|_{L_{x}^{2}}\d \tau.
	\end{aligned}
	\end{equation*}
	Using~\eqref{est:GammaGammahat} and~\eqref{est:Boot21}, we have
	\begin{equation*}
	\mathcal{J}_{9}\lesssim \sum_{\substack{|I_{1}|\le N\\|I_{2}|\le N-6}}\int_{2}^{t}\big\|\widehat{Z}^{I_{1}}\psi\big\|_{L_{x}^{2}}\left\|Z^{I_{2}}\psi\right\|_{L_{x}^{\infty}}\d \tau\lesssim C_{0}^{2}\varepsilon^{2}\int_{2}^{t}\tau^{-\frac{3}{2}}\d \tau \lesssim C_{0}^{2}\varepsilon^{2},
	\end{equation*}
	\begin{equation*}
	\begin{aligned}
	\mathcal{J}_{10}
	&\lesssim \sum_{|I_{1}|\le N-1}\sum_{\substack{|I_{2}|\le N-6\\|I_{3}|\le N-6}}\int_{2}^{t}\left\|Z^{I_{1}}\phi\right\|_{L_{x}^{2}}\left\|(\tau+r){Z}^{I_{2}}\psi\right\|_{L_{x}^{\infty}}\left\|Z^{I_{3}}\psi\right\|_{L_{x}^{\infty}}\d \tau\\
	&\lesssim C_{0}^{3}\varepsilon^{3}\int_{2}^{t}\tau^{-\frac{1}{2}}\tau^{-\frac{3}{2}}\d \tau 
	\lesssim C_{0}^{3}\varepsilon^{3}\int_{2}^{t}\tau^{-2}\d \tau\lesssim C_{0}^{3}\varepsilon^{3},
	\end{aligned}
	\end{equation*}
		\begin{equation*}
	\begin{aligned}
	\mathcal{J}_{11}
	&\lesssim \sum_{|I_{1}|\le N-7}\sum_{\substack{|I_{2}|\le N-6\\|I_{2}|\le N}}\int_{2}^{t}\left\|Z^{I_{1}}\phi\right\|_{L_{x}^{\infty}}\left\|(\tau+r){Z}^{I_{2}}\psi\right\|_{L_{x}^{\infty}}\big\|\widehat{Z}^{I_{3}}\psi\big\|_{L_{x}^{2}}\d \tau\\
	&\lesssim C_{0}^{3}\varepsilon^{3}\int_{2}^{t}\tau^{-1}\tau^{-\frac{1}{2}}\d \tau 
	\lesssim C_{0}^{3}\varepsilon^{3}\int_{2}^{t}\tau^{-\frac{3}{2}}\d \tau\lesssim C_{0}^{3}\varepsilon^{3},
	\end{aligned}
	\end{equation*}
		\begin{equation*}
	\begin{aligned}
	\mathcal{J}_{12}
	&\lesssim\sum_{\substack{|I_{1}|\le N-1\\ |I_{2}|\le N-7}}\sum_{\substack{|I_{3}|\le N-6\\|I_{4}|\le N-6}}\int_{2}^{t} \left\|Z^{I_{1}}\phi\right\|_{L_{x}^{2}} \left\|Z^{I_{2}}\phi\right\|_{L_{x}^{\infty}}\left\|(\tau+r){Z}^{I_{3}}\psi\right\|_{L_{x}^{\infty}}\left\|Z^{I_{4}}\psi\right\|_{L_{x}^{\infty}}\d \tau\\
	&\lesssim C_{0}^{4}\varepsilon^{4}\int_{2}^{t}\tau^{-1}\tau^{-\frac{1}{2}}\tau^{-\frac{3}{2}}\d \tau 
	\lesssim C_{0}^{4}\varepsilon^{4}\int_{2}^{t}\tau^{-3}\d \tau\lesssim C_{0}^{4}\varepsilon^{4},
	\end{aligned}
	\end{equation*}
		\begin{equation*}
	\begin{aligned}
	\mathcal{J}_{13}
	&\lesssim\sum_{\substack{|I_{1}|\le N-7\\ |I_{2}|\le N-7}}\sum_{\substack{|I_{3}|\le N-1\\|I_{4}|\le N-6}}\int_{2}^{t}\left\|Z^{I_{1}}\phi\right\|_{L_{x}^{\infty}}\left\|Z^{I_{2}}\phi\right\|_{L_{x}^{\infty}}\big\|\widehat{Z}^{I_{3}}\psi\big\|_{L_{x}^{2}}\left\|(\tau+r)Z^{I_{4}}\psi\right\|_{L_{x}^{2}}\d \tau\\
	&\lesssim C_{0}^{4}\varepsilon^{4}\int_{2}^{t}\tau^{-1}\tau^{-1}\tau^{-\frac{1}{2}}\d \tau 
	\lesssim C_{0}^{4}\varepsilon^{4}\int_{2}^{t}\tau^{-\frac{5}{2}}\d \tau\lesssim C_{0}^{4}\varepsilon^{4}.
	\end{aligned}
	\end{equation*}
	Combining the above inequalities, we obtain~\eqref{est:N1N2}.
	\end{proof}

Last, we introduce the $L_{t}^{1}L_{x}^{2}$ estimate for the nonlinear terms which will be needed in the energy estimates of wave component $\phi$.
\begin{lemma}
	For all $t\in [2,T^{*}(\psi_{0},\vec{\phi}_{0}))$, we have 
	\begin{equation}\label{est:energyphi}
	\sum_{|I|\le N}\int_{2}^{t}\left\|Z^{I}\left(\psi^{*}\Upsilon\psi\right)\right\|_{L_{x}^{2}}\d \tau\lesssim C_{0}^{2}\varepsilon^{2}.
	\end{equation}
\end{lemma}

\begin{proof}
	From~\eqref{est:Boot21}, $N\ge 13$ and the Leibniz rule, we have 
	\begin{equation*}
	\begin{aligned}
	\sum_{|I|\le N}\int_{2}^{t}\left\|Z^{I}\left(\psi^{*}\Upsilon\psi\right)\right\|_{L_{x}^{2}}\d \tau
	&\lesssim \sum_{\substack{|I_{1}|\le N\\|I_{2}|\le N-7}}\int_{2}^{t}\left\|Z^{I_{1}}\psi\right\|_{L_{x}^{2}}\left\|Z^{I_{2}}\psi\right\|_{L_{x}^{\infty}}\d \tau\\
	&\lesssim C_{0}^{2}\varepsilon^{2}\int_{2}^{t}\tau^{-\frac{3}{2}}\d \tau \lesssim C_{0}^{2}\varepsilon^{2},
	\end{aligned}
	\end{equation*}
	which means~\eqref{est:energyphi}.
	\end{proof}

\subsection{End of the proof of Theorem~\ref{thm:massive}}\label{SS:EndThm2}
We are in a position to complete the proof of Theorem~\ref{thm:massive} and start with the proof of Proposition~\ref{pro:main2}.

\begin{proof}[Proof of Proposition~\ref{pro:main2}]
	For any initial data $(\psi_{0},\vec{\phi}_{0})$ satisfying the smallness condition~\eqref{est:small2}, we consider the corresponding solution $(\psi,\phi)$ of~\eqref{equ:3DDKG2}. From the smallness condition~\eqref{est:small2}, we know that 
	\begin{equation}\label{est:init2}
	\begin{aligned}
	\sum_{|I|\le N}\left(\mathcal{E}(2,\Gamma^{I}\tilde{\phi})^{\frac{1}{2}}+\mathcal{F}(2,\Gamma^{I}\tilde{\phi})^{\frac{1}{2}}\right)&\lesssim \varepsilon,\\
	\sum_{|I|\le N}\left(\mathcal{E}_{D}^{\delta}\big(2,\widehat{Z}^{I}\psi\big)^{\frac{1}{2}}+\mathcal{E}\left(2,Z^{I}\phi\right)^{\frac{1}{2}}\right)&\lesssim \varepsilon.
	\end{aligned}
	\end{equation}
	
		In what follows, we prove Proposition~\ref{pro:main2} by improving all the estimates of $(\psi,\phi)$ in the bootstrap assumption~\eqref{est:Boot21}.
		
		\textbf{Step 1.} Closing the pointwise estimate of $\psi$. Recall that, from~\eqref{equ:3DDKG2} and~\eqref{equ:comm},
		\begin{equation*}
		-\Box Z^{I}\psi+Z^{I}\psi=-Z^{I}\left((\partial_{\mu}\phi)\gamma^{\mu}\gamma^{5}\psi+\phi^{2}\psi\right).
		\end{equation*}
		From~\eqref{est:Lx2psi}, we find
		\begin{equation*}
\sum_{|I|\le N-2}	\max_{2\le \tau\le t}	\tau^{\frac{1}{2}}\left\|Z^{I}\left((\partial_{\mu}\phi)\gamma^{\mu}\gamma^{5}\psi+\phi^{2}\psi\right)\right\|_{L_{x}^{2}}\lesssim C_{0}^{2}\varepsilon^{2}.
		\end{equation*}
		It follows from~\eqref{est:small2} and Corollary~\ref{co:Klein} that 
		\begin{equation*}
		\langle t+r\rangle^{\frac{3}{2}}\sum_{|I|\le N-6}\left|Z^{I}\psi(t,x)\right|\lesssim \varepsilon+ C_{0}^{2}\varepsilon^{2}.
		\end{equation*}
		This inequality strictly improves the pointwise estimates of $\psi$ in the bootstrap assumption~\eqref{est:Boot21} for $C_{0}$ large enough and $\varepsilon$ small enough.
		
		\smallskip
		\textbf{Step 2.} Closing the pointwise estimate of $\phi$. Recall that, in Lemma~\ref{le:Hidden}, we introduce $\tilde{\phi}:=\phi+\frac{1}{4}\psi^{*}\Upsilon\psi$ and obtain
		\begin{equation*}
		-\Box \tilde{\phi}=\frac{1}{2}\left(\mathcal{N}_{1}+\mathcal{N}_{2}+\mathcal{N}_{3}\right).
		\end{equation*}
		We decompose the wave component $\tilde{\phi}$ as,
		\begin{equation*}
		\tilde{\phi}(t,x)=\tilde{\phi}_{free}(t,x)+\tilde{\phi}_{inho}(t,x),
		\end{equation*}
			where $\tilde{\phi}_{free}$ and $\tilde{\phi}_{inho}$ are the solutions for the following 3D linear homogeneous
		or inhomogeneous wave equations,
		\begin{equation*}
		\left\{\begin{aligned}
		-\Box \tilde{\phi}_{free}&=0, \quad \quad \quad \quad \quad \quad \quad \quad \ \mbox{with}\ \ (\tilde{\phi}_{free},\partial_{t}\tilde{\phi}_{free})_{|t=2}=(\tilde{\phi},\pt \tilde{\phi})_{|t=2},\\
		-\Box \tilde{\phi}_{inho}&=\frac{1}{2}\left(\mathcal{N}_{1}+\mathcal{N}_{2}+\mathcal{N}_{3}\right),\quad  \mbox{with}\ \ (\tilde{\phi}_{inho},\partial_{t}\tilde{\phi}_{inho})_{|t=2}=(0,0).
		\end{aligned}\right.
		\end{equation*}
		First, from Lemma~\ref{le:waveenergy} and~\eqref{est:init2}, we have 
		\begin{equation*}
		\sum_{|I|\le N-5}\left\|\Gamma^{I}\tilde{\phi}_{free}\right\|_{L_{x}^{2}}\lesssim \sum_{|I|\le N-6}\left(\mathcal{E}(2,\Gamma^{I}\tilde{\phi})^{\frac{1}{2}}+\mathcal{F}(2,\Gamma^{I}\tilde{\phi})^{\frac{1}{2}}\right)\lesssim \varepsilon.
		\end{equation*}
		It follows from~\eqref{est:Klainer} that
		\begin{equation}\label{est:phiend2}
		\sum_{|I|\le N-7}\big|Z^{I}\tilde{\phi}_{free}\big|\lesssim \varepsilon\langle t+r\rangle^{-1}\langle t-r\rangle^{-\frac{1}{2}}.
		\end{equation}
	
	Second, from~\eqref{est:Qab} and~\eqref{est:Boot21}, we see that 
	\begin{equation*}
\sum_{|I|\le N-7}\left(	\left|Z^{I}\mathcal{N}_{1}\right|+\left|Z^{I}\mathcal{N}_{2}\right|+\left|Z^{I}\mathcal{N}_{3}\right|\right)\lesssim C_{0}^{2}\varepsilon^{2}(t+r)^{-4}.
	\end{equation*} 
	Therefore, from Proposition~\ref{pro:Linwave}, Lemma~\ref{le:waveenergy},~\eqref{est:small2} and~\eqref{est:Klainer}, we have 
	\begin{equation}\label{est:phiend22}
	\sum_{|I|\le N-7}\big|Z^{I}\tilde{\phi}_{inho}\big|\lesssim 
	\varepsilon\langle t+r\rangle^{-1}\langle t-r\rangle^{-\frac{1}{2}}+
	C_{0}^{2}\varepsilon^{2}\langle t+r\rangle^{-1}\langle t-r\rangle^{-\frac{1}{2}}.
	\end{equation}
	Combining~\eqref{est:Boot21},~\eqref{est:phiend2},~\eqref{est:phiend22} with the definition of $\tilde{\phi}$, we obtain
	\begin{equation*}
	\begin{aligned}
	\sum_{|I|\le N-7}\left|Z^{I}\phi\right|
	&\lesssim \sum_{|I|\le N-7}\left(\big|Z^{I}\tilde{\phi}_{free}\big|+\big|Z^{I}\tilde{\phi}_{inho}\big|\right)+\sum_{\substack{|I_{1}|\le N-7\\ |I_{2}|\le N-7}}\left|Z^{I_{1}}\psi\right|\left|Z^{I_{2}}\psi\right|\\
	&\lesssim \left(\varepsilon+C_{0}^{2}\varepsilon^{2}\right)\langle t+r\rangle^{-1}\langle t-r\rangle^{-\frac{1}{2}}+C_{0}^{2}\varepsilon^{2}\langle t+r\rangle^{-3}.
	\end{aligned}
	\end{equation*}
	This inequality strictly improves the pointwise estimates of $\phi$ in the bootstrap assumption~\eqref{est:Boot21} for $C_{0}$ large enough and $\varepsilon$ small enough.
	
	\smallskip
	\textbf{Step 3.} Bound on the weighted $L_{x}^{2}$ and energy norms of $\psi$. First, from ~\eqref{equ:3DDKG2vec}, \eqref{est:Lt1Lx2psi}, \eqref{est:init2} and Lemma~\ref{le:DiracGhotst}, we have 
	\begin{equation*}
	\begin{aligned}
	\sum_{|I|\le N}\mathcal{E}_{D}^{\delta}(t,\widehat{Z}^{I}\psi)
	&\lesssim \sum_{|I|\le N}\int_{2}^{t}\int_{\R^{3}}\big|\big(\widehat{Z}^{I}\psi\big)\gamma^{0}\widehat{Z}^{I}\big(i\phi \gamma^{5}\psi\big)\big| \d x \d \tau\\
	&+\sum_{|I|\le N}\mathcal{E}_{D}^{\delta}(2,\widehat{Z}^{I}\psi)\lesssim \varepsilon^{2}+C_{0}^{3}\varepsilon^{3},
	\end{aligned}
	\end{equation*}
	which implies
	\begin{equation}\label{est:Diracimpro}
	\sum_{|I|\le N}\mathcal{E}_{D}^{\delta}(t,\widehat{Z}^{I}\psi)^{\frac{1}{2}}\lesssim \left(\varepsilon^{2}+C_{0}^{3}\varepsilon^{3}\right)^{\frac{1}{2}}\lesssim \varepsilon+C_{0}^{\frac{3}{2}}\varepsilon^{\frac{3}{2}}.
	\end{equation}

	Second, we recall that, for any $I\in \mathbb{N}^{10}$ with $|I|\le 10$,
		\begin{equation*}
	-\Box Z^{I}\psi+Z^{I}\psi=-Z^{I}\left((\partial_{\mu}\phi)\gamma^{\mu}\gamma^{5}\psi+\phi^{2}\psi\right).
	\end{equation*}
	Therefore, from~\eqref{est:GammaGammahat},~\eqref{est:Boot21},~\eqref{est:Lx2psi} and Lemma~\ref{le:decayKG}, we obtain
	\begin{equation}\label{est:DiracL2impro}
	\begin{aligned}
	&\sum_{|I|\le N-2}\left\|\frac{\langle t+r\rangle}{\langle t-r\rangle}Z^{I}\psi\right\|_{L_{x}^{2}}\\
	&\lesssim \sum_{|I|\le N}\big\|\widehat{Z}^{I}\psi\big\|_{L_{x}^{2}}+\sum_{|I|\le N-2} \left\|(t+r)Z^{I}(\psi \partial \phi)\right\|_{L^{2}_{x}}\\
	&+\sum_{|I|\le N-2} \left\|(t+r)Z^{I}(\phi^{2}\psi )\right\|_{L^{2}_{x}}\lesssim 
	\varepsilon+C_{0}^{\frac{3}{2}}\varepsilon^{\frac{3}{2}}+C_{0}^{2}\varepsilon^{2}+C_{0}^{3}\varepsilon^{3}.
	\end{aligned}
	\end{equation}
	
	The inequalities~\eqref{est:Diracimpro}-\eqref{est:DiracL2impro} strictly improve the weighted $L_{x}^{2}$ and energy estimates of $\psi$ in the bootstrap assumption~\eqref{est:Boot21} for $C_{0}$ large enough and $\varepsilon$ small enough.

	\smallskip
	\textbf{Step 4.} Bound on the $L_{x}^{2}$ and energy norms of $\phi$. First, from~\eqref{equ:3DDKG2vec},~\eqref{est:energyphi},~\eqref{est:init2} and Lemma~\ref{le:waveenergy}, we see that 
	\begin{equation}\label{est:energyphiimpove}
	\begin{aligned}
	\sum_{|I|\le N}\mathcal{E}(t,Z^{I}\phi)^{\frac{1}{2}}
	&\lesssim \sum_{|I|\le N}\mathcal{E}(2,Z^{I}\phi)^{\frac{1}{2}}\\
	&+\sum_{|I|\le N}\int_{2}^{t}\left\|Z^{I}\left(\psi^{*}\Upsilon\psi\right)\right\|_{L_{x}^{2}}\d \tau\lesssim \varepsilon+C_{0}^{2}\varepsilon^{2}.
	\end{aligned}
	\end{equation}

	Second, from~\eqref{equ:comm} and~\eqref{equ:Hiddenphi}, for any $I\in \mathbb{N}^{10}$ with $|I|\le N-1$
	\begin{equation*}
	-\Box Z^{I}\tilde{\phi}=\frac{1}{2}Z^{I}\mathcal{N}_{1}+\frac{1}{2}Z^{I}\mathcal{N}_{2}+\frac{1}{2}Z^{I}\mathcal{N}_{3}.
	\end{equation*}
	Based on the above identity,~\eqref{est:N1N2} and Lemma~\ref{le:waveenergy}, we have 
	\begin{equation*}
	\begin{aligned}
	&\sum_{|I|\le N-1}\left(\mathcal{E}(t,Z^{I}\tilde{\phi})^{\frac{1}{2}}+\mathcal{F}(t,Z^{I}\tilde{\phi})^{\frac{1}{2}}\right)\\
	&\lesssim 	\sum_{|I|\le N-1}\left(\mathcal{E}(2,Z^{I}\tilde{\phi})^{\frac{1}{2}}+\mathcal{F}(2,Z^{I}\tilde{\phi})^{\frac{1}{2}}\right)+\sum_{|I|\le N-1}\int_{2}^{t}\left\|(\tau+r)Z^{I}\mathcal{N}_{1}\right\|_{L_{x}^{2}}\d \tau\\
	&+\sum_{|I|\le N-1}\int_{2}^{t}\left(\left\|(\tau+r)Z^{I}\mathcal{N}_{2}\right\|_{L_{x}^{2}}+\left\|(\tau+r)Z^{I}\mathcal{N}_{3}\right\|_{L_{x}^{2}}\right)\d \tau \lesssim \varepsilon+C_{0}^{2}\varepsilon^{2}+C_{0}^{4}\varepsilon^{4},
	\end{aligned}
	\end{equation*}
	which implies
	\begin{equation*}
	\sum_{|I|\le N}\big\|Z^{I}\tilde{\phi}\big\|_{L_{x}^{2}}\lesssim \sum_{|I|\le N-1}\left(\mathcal{E}(t,Z^{I}\tilde{\phi})^{\frac{1}{2}}+\mathcal{F}(t,Z^{I}\tilde{\phi})^{\frac{1}{2}}\right)\lesssim \varepsilon+C_{0}^{2}\varepsilon^{2}+C_{0}^{4}\varepsilon^{4}.
	\end{equation*}
	Therefore, using again~\eqref{est:Boot21} and the definition of $\tilde{\phi}$, we obtain
	\begin{equation}\label{est:L2phiimprove}
	\begin{aligned}
	\sum_{|I|\le N}\big\|Z^{I}{\phi}\big\|_{L_{x}^{2}}
	\lesssim &\sum_{\substack{|I_{1}|\le N\\|I_{2}|\le N-6}}\left\|Z^{I_{1}}\psi\right\|_{L_{x}^{2}}\left\|Z^{I_{2}}\psi\right\|_{L_{x}^{\infty}}\\
	&+\sum_{|I|\le N}\big\|Z^{I}\tilde{\phi}\big\|_{L_{x}^{2}}\lesssim \varepsilon+C_{0}^{2}\varepsilon^{2}+C_{0}^{4}\varepsilon^{4}.
	\end{aligned}
	\end{equation}
	The inequalities~\eqref{est:energyphiimpove}-\eqref{est:L2phiimprove} strictly improve the $L_{x}^{2}$ and energy estimates of $\phi$ in the bootstrap assumption~\eqref{est:Boot21} for $C_{0}$ large enough and $\varepsilon$ small enough.

		\smallskip
	At this point, we have strictly improved all the pointwise, $L_{x}^{2}$ and energy estimates of $(\psi,\phi)$ in the bootstrap assumption~\eqref{est:Boot21}.  In conclusion, for all initial data $(\psi_{0},\vec{\phi}_{0})$ satisfying~\eqref{est:small2}, we show that $T^{*}(\psi_{0},\vec{\phi}_{0})=\infty$ and thus the proof of Proposition~\ref{pro:main2} is complete.

\end{proof}

Now, we complete the proof of Theorem~\ref{thm:massive} by using Proposition~\ref{pro:main2}.

\begin{proof}[End of the proof of Theorem~\ref{thm:massive}]
	Notice that, from~\eqref{est:Boot21}, the global
	existence of $(\psi,\phi)$, the pointwise decay~\eqref{est:theorem2point} and the energy bound~\eqref{est:theorem2} are consequences of Proposition~\ref{pro:main2}.
	Our task now is to prove that the solution scatters to a free solution in the
	energy space $H^{N-2}\times \mathcal{X}^{N}$.
	
	\smallskip
	For the scattering of the Dirac spinor $\psi$, we first recall a few facts about the linear theory. Consider the Cauchy problem for the 3D linear Dirac equation, 
	\begin{equation}\label{equ:Diracsca}
	-i\gamma^{\mu}\partial_{\mu}u+u=G\quad \mbox{with}\quad u_{|t=2}=u_{0}.
	\end{equation}
	The Fourier analysis gives us a solution of~\eqref{equ:Diracsca},
	\begin{equation}\label{equ:DuhamDirac}
	u(t)=\mathcal{S}_{D}(t-2)u_{0}
	+i\int_{2}^{t}\mathcal{S}_{D}(t-\tau)G(\tau)\d \tau,
	\end{equation}
	where $\mathcal{S}_{D}(t)$ can be written as 
	\begin{equation*}
	\mathcal{S}_{D}(t)=
	e^{it(i\gamma^{0}\gamma^{a}\partial_{a}+\gamma^{0})},\quad \forall t\in (-\infty,\infty).
	\end{equation*}
	Here, the operators mean that 
	\begin{equation*}
	\left(\mathcal{S}_{D}(t)f\right)(x)=\mathcal{F}^{-1}\left(e^{it\left(-\gamma^{0}\gamma^{a}\xi_{a}+\gamma^{0}\right)}\hat{f}(\xi)\right)(x).
	\end{equation*}
	 Recall that, for $K\in \mathbb{N}$, the 3D Dirac flow $\mathcal{S}_{D}(t)$ are isometric groups on $H^{K}$. More precisely, for $K\in \mathbb{N}$, we have the following identities,
	\begin{equation}\label{est:Sd}
	\left\|\mathcal{S}_{D}(t)\right\|_{{H}^{K}}=1,\quad \forall t\in (-\infty,\infty).
	\end{equation}
	Recall also that, we have the following basic properties of $\mathcal{S}_{D}$,
	\begin{equation}\label{equ:Stsd}
	\mathcal{S}_{D}(t+\tau)=\mathcal{S}_{D}(t)\circ\mathcal{S}_{D}(\tau),\quad \forall(t,\tau)\in (-\infty,\infty)^{2}.
	\end{equation}
	
	\smallskip
	Now, we show that the Dirac spinor $\psi$ scatters to a free solution. In actually, from~\eqref{est:pGammaGammahat},~\eqref{est:Boot21} and $T_{*}(\psi_{0},\vec{\phi}_{0})=\infty$, we have
	\begin{equation*}
	\begin{aligned}
	&\sum_{|I|\le N-2}\int_{2}^{\infty}\big\|\widehat{Z}^{I}(i\phi\gamma^{5}\psi)\big\|_{L_{x}^{2}}\d \tau\\
	&\lesssim \sum_{\substack{|I_{1}|\le N-7\\|I_{2}|\le N-2}}\int_{2}^{\infty}\tau^{-1}\left\|\langle \tau -r \rangle Z^{I_{1}}\phi\right\|_{L_{x}^{\infty}}\left\|\frac{\langle\tau+r\rangle}{\langle \tau-r\rangle}Z^{I_{2}}\psi\right\|_{L_{x}^{2}}\d \tau\\
	&+\sum_{\substack{|I_{1}|\le N-2\\|I_{2}|\le N-6}}\int_{2}^{\infty}\left\|Z^{I_{1}}\phi\right\|_{L_{x}^{2}}\left\|Z^{I_{2}}\psi\right\|_{L_{x}^{\infty}}\d \tau\lesssim C_{0}^{2}\varepsilon^{2},
	\end{aligned}
	\end{equation*}
	which implies
	\begin{equation}\label{est:psiHN-2}
\int_{2}^{\infty}\big\|(i\phi\gamma^{5}\psi)\big\|_{H^{N-2}}\d \tau\lesssim\sum_{|I|\le N-2}\int_{2}^{\infty}\big\|\widehat{Z}^{I}(i\phi\gamma^{5}\psi)\big\|_{L_{x}^{2}}\d \tau\lesssim C_{0}^{2}\varepsilon^{2}.
	\end{equation}
	Based on the above inequality and~\eqref{est:Sd}, we have 
	\begin{equation*}
	\int_{2}^{\infty}\left\|\mathcal{S}_{D}(2-\tau)(i\phi\gamma^{5}\psi)\right\|_{H^{N-2}}\d \tau \lesssim \int_{2}^{\infty}\left\|i\phi\gamma^{5}\psi\right\|_{H^{N-2}}\d \tau\lesssim C\varepsilon,
	\end{equation*}
	which implies
	\begin{equation*}
	\psi_{0h}:=\int_{2}^{\infty}\mathcal{S}_{D}(2-\tau)(i\phi\gamma^{5}\psi)\d \tau \in H^{N-2}.
	\end{equation*}
	Therefore,
	\begin{equation*}
	\psi_{0\ell}:=\psi_{0}+\psi_{0h}\in H^{N-2}\Rightarrow 
	\psi_{\ell}(t)=\mathcal{S}_{D}(t-2)\psi_{0\ell}\in H^{N-2}.
	\end{equation*}
	Last, from~\eqref{est:Sd},~\eqref{equ:Stsd} and~\eqref{est:psiHN-2}, we obtain
	\begin{equation*}
	\lim_{t\to \infty}\left\|\psi(t)-\psi_{\ell}(t)\right\|_{H^{N-2}}\le \lim_{t\to \infty}\int_{t}^{\infty}\left\|i\phi\gamma^{5}\psi\right\|_{H^{N-2}}\d \tau=0,
	\end{equation*}
	which means that the Dirac spinor $\psi$ scatters to a free solution in $H^{N-2}$. 
	
	For the scattering of the wave component $\phi$, we use the similar argument as in the end of the proof of Theorem~\ref{thm:massless}. More precisely, from~\eqref{est:energyphi}, we see that 
	\begin{equation*}
	\sum_{|I|\le N}\int_{2}^{t}\left\|\psi^{*}\Upsilon\psi\right\|_{H^{N}}\d \tau\lesssim
		\sum_{|I|\le N}	\int_{2}^{t}\left\|Z^{I}\left(\psi^{*}\Upsilon\psi\right)\right\|_{L_{x}^{2}}\d \tau\lesssim C_{0}^{2}\varepsilon^{2}.
	\end{equation*}
	Therefore, using the linear theory of the 3D wave equation, we conclude that there exists $(\phi_{0\ell},\phi_{1\ell})\in \mathcal{X}^{N}$ such that 
	\begin{equation*}
	\lim_{t\to \infty}\left\|(\phi,\pt \phi)-(\phi_{\ell},\pt \phi_{\ell})\right\|_{\mathcal{X}^{N}}=0,
\end{equation*}
where $\phi_{\ell}$ is the solution to the 3D linear wave equation with the initial data $(\phi_{0\ell},\phi_{1\ell})$. In conclusion, the solution $(\psi,\phi)$ has the required properties and so the proof of Theorem~\ref{thm:massive} is complete.
\end{proof}

\end{document}